\documentclass[11pt,a4paper]{article}
\input glyphtounicode
\usepackage[T1]{fontenc}
\usepackage[utf8]{inputenc}
\usepackage{lmodern}
\usepackage{amsmath,amssymb,amsthm,mathrsfs}
\usepackage[margin=27mm,headheight=14pt]{geometry}
\usepackage{booktabs,longtable,array,enumitem,xcolor}
\usepackage{fancyhdr}
\usepackage[unicode=true,colorlinks=true,linkcolor=blue!50!black,citecolor=blue!50!black,urlcolor=blue!50!black]{hyperref}
\usepackage{xurl}
\hypersetup{pdftitle={Center-uniform uncertainty for the Dunkl transform},pdfauthor={Xingyu Zhao and Zhiwen Duan}}
\numberwithin{equation}{section}
\newtheorem{theorem}{Theorem}[section]
\newtheorem{lemma}[theorem]{Lemma}
\newtheorem{proposition}[theorem]{Proposition}
\newtheorem{corollary}[theorem]{Corollary}
\theoremstyle{definition}
\newtheorem{definition}[theorem]{Definition}
\theoremstyle{remark}
\newtheorem{remark}[theorem]{Remark}
\DeclareMathOperator{\supp}{supp}
\DeclareMathOperator{\dist}{dist}
\DeclareMathOperator{\conv}{conv}
\newcommand{\Fk}{\mathcal F_k}
\newcommand{\R}{\mathbb{R}}
\newcommand{\C}{\mathbb{C}}

\newcommand{\doi}[1]{\href{https://doi.org/#1}{\nolinkurl{doi:#1}}}
\hypersetup{bookmarksdepth=3}
\allowdisplaybreaks[2]
\setlist{itemsep=2pt,topsep=5pt}
\title{Center-uniform uncertainty for the Dunkl transform and  a fractal uncertainty principle}
\author{Xingyu Zhao\thanks{Corresponding author.}\quad Zhiwen Duan\\[3pt]
\small School of Mathematics and Statistics,\\
\small Huazhong University of Science and Technology,\\
\small Wuhan 430074, P.R. China\\[3pt]
\small\texttt{zhaoxingyu@hust.edu.cn}\quad\texttt{duanzhw@hust.edu.cn}}

\date{}
\begin{document}
\maketitle
\begin{abstract}
For a fixed finite reflection group and nonnegative multiplicity, we prove
that quasianalytic Dunkl spectral decay implies observation on arbitrary
thick sets, uniformly in the frequency center. This yields a
Logvinenko--Sereda theorem for translated spectral balls. We also prove
necessity of the logarithmic integral condition at every spectral center.
We establish a fractal uncertainty principle for
spatial sets porous on balls and spectral sets porous on lines, with
concentration bound $Ch^\beta$, $\beta>0$.
 
\end{abstract}

\section{Introduction}\label{sec:introduction}
\subsection{Notation}\label{subsec:notation}
Write $B(x,r)=\{y:|y-x|<r\}$ and $Q_L=[-L/2,L/2]^d$.

We use the unitary Fourier transform
$\widehat f(\xi)=(2\pi)^{-d/2}\int e^{-i\langle x,\xi\rangle}f(x)\,dx$.
For a fixed root system $R$, reflection group $G$, and nonnegative
$G$-invariant multiplicity $k$, let
\[
 d\mu_k(x)=\prod_{\alpha\in R_+}|\langle\alpha,x\rangle|^{2k_\alpha}\,dx,
 \qquad L_k^2=L^2(\R^d,\mu_k),\qquad
 \|f\|_{2,k}=\|f\|_{L_k^2}.
\]
 
Here $R_+$ is a positive subsystem of $R$ (see Section 2 for details),
and $\Fk$ denotes the unitary
Dunkl transform~\cite{MFE1,MR1}. Write $gE=\{gx:x\in E\}$. A set $E$ is invariant if $gE=E$ for
all $g\in G$. All supports are essential supports.

Throughout, constants may depend on the fixed root system $R$ and
multiplicity $k$; this dependence is suppressed unless stated otherwise.

\subsection{Uncertainty principles and main results}
\begin{definition}\label{def:thick}
Let $0<\gamma\le1$ and $L>0$. A measurable set $E\subset\R^d$ is
$(\gamma,L)$-thick with respect to $\mu$ if
\[
 \mu(E\cap(x+Q_L))\ge\gamma\,\mu(x+Q_L)\qquad(x\in\R^d).
\]
We use this definition for Lebesgue measure and for $\mu_k$.
\end{definition}
The Logvinenko--Sereda theorem~\cite{LVN} (see also~\cite{FS}) gives the following estimate for band-limited functions:
for a Lebesgue $(\gamma,L)$-thick set $E$,
\begin{equation}\label{eq:intro-fourier-logvinenko-sereda}
 \supp\widehat f\subset B(0,N)
 \quad\Longrightarrow\quad
 \|f\|_{L^2(\R^d)}\le C(\gamma,LN,d)\|f\|_{L^2(E)}.
\end{equation}
Kovrijkine~\cite{OK} obtained quantitative bounds for this estimate.
The relation between spatial regularity and Fourier decay suggests replacing
compact spectral support by a weighted decay condition. Jaye and
Mitkovski~\cite{BJ} established observation under
\eqref{eq:intro-fourier-decay} below and characterized the admissible weights
by the logarithmic integral condition, under the stated regularity assumptions.

\begin{equation}\label{eq:intro-fourier-decay}
 \|W(|\cdot|)\widehat f\|_2\le C_W\|f\|_2,
\end{equation}
where $C_W$ is fixed and $W$ is a unique continuation weight, defined as follows.

\begin{definition}\label{def:UCweight}
A \emph{unique continuation weight} is a nondecreasing, lower
semicontinuous function $W:[0,\infty)\to[1,\infty]$ with $W(0)=1$
satisfying the following conditions:
\begin{enumerate}[label=(\roman*),nosep]
\item $s\mapsto\log W(e^s)$ is convex on $[0,\infty)$.
\item For each integer $n\ge0$ there is $c_n>0$ such that
      $W(t)\ge c_n t^n$ for all $t\ge0$.
\item $\displaystyle\int_0^\infty\frac{\log W(t)}{1+t^2}\,dt=\infty$.
\end{enumerate}
\end{definition}
This is the weight class of Cohen~\cite[Definition~3.1]{AC},
following the quantitative uniqueness theorem of Jaye and
Mitkovski~\cite[Theorem~1.3]{BJ}.

The quasianalytic condition on $W$ links spectral moments to
propagation of smallness and hence to an observation inequality
$\|f\|_2\le C(W,C_W,\gamma,L,d)\|f\|_{L^2(E)}$.

In addition, if we take such a weight $W$, 
\begin{equation}\label{eq:example}
  W(t)=
\begin{cases}
1 & \text{for } |t|\le N,\\
+\infty & \text{for } |t|>N,
\end{cases}  
\end{equation}

then this reduces to \eqref{eq:intro-fourier-logvinenko-sereda}; we refer the reader to \cite{BJ} for further details.

The Fourier uncertainty principle limits simultaneous localization
of a function and its Fourier transform. In this sense,
\eqref{eq:intro-fourier-logvinenko-sereda} states that the full energy of a
band-limited function is controlled by its energy on a Lebesgue
$(\gamma,L)$-thick set $E$. Related uncertainty estimates concern
fractal sets, whose geometry is quantified here by porosity.

\begin{definition}\label{def:porosity}
Let $\nu>0$ and $0<\alpha_0\le\alpha_1$.
A set $X\subset\R^d$ is \emph{$\nu$-porous on balls} from scales
$\alpha_0$ to $\alpha_1$ if for each ball $B$ of diameter
$r\in[\alpha_0,\alpha_1]$ there exists $x\in B$ with
$B(x,\nu r)\cap X=\varnothing$.
It is \emph{$\nu$-porous on lines} over the same scales if each
line segment $\tau$ of length $r\in[\alpha_0,\alpha_1]$ contains
a point $x$ with $B(x,\nu r)\cap X=\varnothing$.
\end{definition}

One obtains an uncertainty principle of the following form: for $0<h<1$,
\begin{equation}\label{eq:intro-fourier-fractal-uncertainty}
 \supp\widehat f\subset Y
 \quad\Longrightarrow\quad
 \|\mathbf1_Xf\|_2\le Ch^\beta\|f\|_2,
 \qquad \beta>0,
\end{equation}
when $X\subset[-1,1]^d$ is porous on balls from scales $h$ to $1$ and
$Y\subset[-h^{-1},h^{-1}]^d$ is porous on lines from scales $1$ to $h^{-1}$.

Dyatlov and Zahl~\cite{DZ} introduced a fractal uncertainty principle
in their study of spectral gaps and additive energy. Bourgain and
Dyatlov~\cite{BD} established a power gain for one-dimensional
regular fractal sets. Dyatlov and Jin~\cite{DJ} developed an
alternative approach based on Dolgopyat's method, and Backus, Leng,
and Tao~\cite{BLT} obtained higher-dimensional results under suitable
geometric hypotheses. Cohen~\cite{Cohen2025} proved the
estimate \eqref{eq:intro-fourier-fractal-uncertainty} using a
higher-dimensional Beurling--Malliavin multiplier theorem; his
thesis~\cite{AC} develops the damping and transference framework
used below. In this approach, damping converts porous spectral
support into weighted spectral decay, observation forces mass into
spatial holes, and iteration produces $h^\beta$.

These uncertainty principles also have applications to evolution
and spectral problems. Quantitative Logvinenko--Sereda estimates are used to prove
null-controllability of the heat equation from thick sets
\cite{EV}, while fractal uncertainty principles yield essential
spectral gaps for convex co-compact hyperbolic surfaces~\cite{BD}.

Dunkl theory originally arose from a concrete mathematical problem: how to extend the classical theory of spherical harmonics and orthogonal polynomials to weighted settings with the symmetry of a finite reflection group. Against this backdrop, Charles F. Dunkl introduced in 1989 a family of differential-difference operators now bearing his name and proved that these operators commute pairwise. His work on reflection-invariant spherical polynomials, differential-difference operators, integral kernels, and Hankel transforms
\cite{CFD0,CFD1,CFD2,CFD3}, together with de Jeu's transform theory
\cite{MFE1}, provides a reflection-group analogue of Fourier analysis.

Specifically, at $k=0$ the Dunkl transform is the ordinary Fourier transform. On radial functions it reduces to a Hankel transform of order $d_k/2-1$, where $d_k$ is the homogeneous dimension in \eqref{eq:weight} (see~\cite[Section 2]{ST1}).
We use the basic framework surveyed in~\cite{MR1}.

Several classical uncertainty principles have been extended to this setting.
R\"osler~\cite{RU} proved a Heisenberg--Weyl inequality,
and Shimeno~\cite{SH} studied Heisenberg and Hardy principles.
Kawazoe and Mejjaoli~\cite{KM} obtained versions of Cowling--Price,
Beurling, and Donoho--Stark uncertainty principles.
Ghobber and Jaming~\cite{GJ} proved local uncertainty and
Benedicks--Amrein--Berthier results for integral transforms including
the Dunkl transform. They also proved a Fourier--Bessel
Logvinenko--Sereda theorem~\cite{SGP}. Ghobber~\cite{GS} established a Shapiro-type
principle for orthonormal bases and quantitative estimates for
orthonormal sequences. 

Recently, Zhao, Wei, and Duan~\cite{S4} obtained an origin-centered Logvinenko--Sereda inequality while studying observability inequalities for the Dunkl heat equation. In the Fourier setting, the modulation $e^{i\langle \eta,x\rangle}$ allows one to shift the spectral support of a band-limited function $f$. That is, the condition $\supp\widehat f\subset B(0,N)$ in \eqref{eq:intro-fourier-logvinenko-sereda} can be naturally replaced by $\supp\widehat f\subset B(\eta,N)$, and the resulting constant $C(\gamma,LN,d)$ is independent of $\eta$. 

However,  in the usual Dunkl setting, this uniformity requires an argument adapted to the reflection structure. Ordinary modulation does not in general translate Dunkl
spectral support. Reflection terms couple different points of an
orbit, and the measure degenerates on reflection hyperplanes. 
This issue constitutes a general difficulty within Dunkl theory. For instance, results concerning the Paley–Wiener theorem (see \cite{MFE2}) and those on transformations of translated supports (see \cite{JDA2}) both impose additional symmetry requirements on the domains or functions involved. All these difficulties point to a fundamental issue: our current understanding of the asymptotic behavior of the Dunkl kernel remains limited.
In particular, when studying the fractal uncertainty principle, we in~\cite{ZWD} used the asymptotics of Bessel functions to establish the fractal uncertainty principle for the Hankel transform. However, the same strategy fails for the fractal uncertainty principle associated with the Dunkl transform.

In this paper, we address the difficulties outlined above and establish two main results. The first extends the Logvinenko–Sereda inequality for the Dunkl transform from~\cite{S4} to a version independent of the location of the spectral support.

\begin{theorem}\label{T91}
Fix $d,W,C_W,\gamma,L$, where $0<\gamma\le1$, $L>0$, $W$ satisfies Definition~\ref{def:UCweight} and $C_W>0$. There exists
\[
C = C(W,C_W,\gamma,L,d)<\infty,
\]
such that for any $\eta\in\mathbb{R}^d$, any measurable set $E$ that is $(\gamma,L)$-thick for $\mu_k$, and any $f\in L_k^2$, we have
\begin{equation}\label{E91}
\bigl\|W(|\cdot-\eta|)\Fk f\bigr\|_{2,k}\le C_W\|f\|_{2,k}
\implies
\|f\|_{2,k}\le C\bigl\|\mathbf 1_E f\bigr\|_{2,k}.
\end{equation}
When $C_W<1$ the spectral condition forces $f=0$, so below we assume $C_W\ge1$. The constant does not depend on the frequency center $\eta$, the specific location of the observation set, or the function $f$. 
\end{theorem}

Taking $W(t)$ as \eqref{eq:example} gives the center-uniform Logvinenko--Sereda theorem
\begin{equation}\label{eq:intro-dunkl-logvinenko-sereda}
 \supp\Fk f\subset\overline{B(\eta,N)}
 \quad\Longrightarrow\quad
 \|f\|_{2,k}\le C(\gamma,LN,d)
                         \|\mathbf1_Ef\|_{2,k}.
\end{equation}

Building upon this, our second main result adopts the strategy from~\cite{BJ} and related works to derive the following fractal uncertainty principle for the Dunkl transform.

\begin{theorem}\label{thm:fractal-uncertainty}
Let $d\ge1$ be an integer, $X,Y$ are measurable sets, let $R\subset\R^d$ be a fixed finite reduced
root system with reflection group $G$, and let $k:R\to[0,\infty)$ be a
fixed $G$-invariant multiplicity. Let $\nu>0$, $0<h<1$, and assume that
\begin{itemize}
\item $\mathbf X\subset[-1,1]^d$ is $\nu$-porous on balls from scales
      $h$ to $1$;
\item $\mathbf Y\subset[-h^{-1},h^{-1}]^d$ is $\nu$-porous on lines
      from scales $1$ to $h^{-1}$.
\end{itemize}
Then there exist $\beta=\beta(\nu,d)>0$ and $C=C(\nu,d)>0$ such that
for all $f\in L_k^2$,
\begin{equation}\label{eq:fractal-uncertainty}
 \supp\Fk f\subset\mathbf Y\quad\Longrightarrow\quad
 \|\mathbf1_{\mathbf X}f\|_{2,k}\le Ch^\beta\|f\|_{2,k}.
\end{equation}
The constants are independent of $h$, $\mathbf X$, $\mathbf Y$, and $f$.
\end{theorem}

\subsection{Main contributions}
We relate quasianalytic Dunkl spectral decay to thick-set observation
uniformly in the spectral center, obtaining the Logvinenko--Sereda
theorem as a special case. Theorem~\ref{thm:necessity} also shows that
the logarithmic integral condition is necessary: when it fails, a common
thick set and a common bound $C_W$ give counterexamples at every spectral
center. Uniform non-concentration near reflection
hyperplanes and quasianalytic propagation for individual orbit components
provide the required reflection-sensitive estimates. Combined with damping and orbit-local
support, this weighted theorem yields the fractal uncertainty principle.

\subsection{Proof strategy}
The weighted observation theorem follows the quasianalytic and good-cube
method of Jaye and Mitkovski~\cite{BJ}, with two adaptations to reflections.
First, a stratified argument excludes concentration near reflection
hyperplanes uniformly in the spectral center. Second, conjugating orbit
values by the square-root weight and scalar phases gives a first-order
system with uniformly bounded coefficients, related to~\cite{MR4}.
In weak limits, oscillatory coefficient terms either disappear or retain
analytic coefficients. Scalar quasianalyticity~\cite{FNM,BJ} of the
limiting components then yields uniform observation of a component
carrying maximal mass. As in the enlarged-set argument of
\cite[proof of Lemma 5.7]{S4}, complete observation cubes and this choice
of component allow arbitrary thick sets. Summation over good cores
gives the estimate. The cutoff weight in \eqref{eq:example} yields
the Logvinenko--Sereda theorem. For necessity, the Paley--Wiener
construction of~\cite[Section~4]{BJ}, together with uniform comparisons
of the translated weighted norms, gives compactly supported
counterexamples at every spectral center.

For the fractal uncertainty principle, Cohen's multiplier theorem
\cite{Cohen2025,AC} and spherical Paley--Wiener theory~\cite{MFE2}
produce uniformly centered damping multipliers. Orbit-local support
\cite{JDA2} and the transference method of~\cite{BJ,AC} give observation
on holes. Positive radial translations~\cite{MR3,ST1} supply cutoffs,
and iteration over scales yields the power gain. Auxiliary orbit
symmetrization is retained for convolution support and these cutoffs.

\subsection{Organization of the paper}
Section~\ref{sec:preliminaries} recalls Dunkl analysis, including thickness,
translation, convolution and the spherical Paley--Wiener theorem. It then
collects the porosity estimates and the Beurling--Malliavin multiplier
theorem, followed by the spectral moment estimates associated with the weights of Definition~\ref{def:UCweight}.
Section~\ref{sec:logvinenko-sereda-main} proves center-uniform observability from
Dunkl spectral decay, obtains the Logvinenko--Sereda theorem for
arbitrary spectral centers, and proves necessity of the logarithmic
integral condition in Section~\ref{subsec:necessity}.
Section~\ref{sec:fractal-uncertainty} constructs the damping multipliers,
proves observation on porous holes, and combines cutoffs with multiscale
iteration to prove Theorem~\ref{thm:fractal-uncertainty}.

\section{Preliminaries}\label{sec:preliminaries}

\subsection{Dunkl analysis}\label{subsec:Dunkl-notation}
Fix a finite reduced root system $R\subset\R^d$, a positive subsystem $R_+$,
and a multiplicity $k:R\to[0,\infty)$ invariant under the reflection group
$G=\langle\sigma_\alpha:\alpha\in R\rangle$. Write $k_\alpha=k(\alpha)$.
We use the weight, measure and homogeneous dimension of Dunkl analysis
as in~\cite[Section 2]{ST1} and~\cite[Section 2]{JDA2,DHtranslations}:
\begin{equation}\label{eq:weight}
 \begin{gathered}
 \sigma_\alpha x=x-2\frac{\langle\alpha,x\rangle}{|\alpha|^2}\alpha,
 \qquad w_k(x)=\prod_{\alpha\in R_+}|\langle\alpha,x\rangle|^{2k_\alpha},
 \qquad d\mu_k(x)=w_k(x)\,dx,\\
 \gamma_k=\sum_{\alpha\in R_+}k_\alpha,\qquad d_k=d+2\gamma_k .
 \end{gathered}
\end{equation}
Factors with $k_\alpha=0$ are omitted. We write $\|\cdot\|_{p,k}$ for the
norm in $L_k^p=L^p(\R^d,\mu_k)$ and $|A|$ for Lebesgue measure.
Here $B_r=B(0,r)$, $Q_L=[-L/2,L/2]^d$, and
$\langle x\rangle=(1+|x|^2)^{1/2}$. The measure $\mu_k$ is homogeneous
of degree $d_k$.

For $v\in\R^d$, the Dunkl operator is
\begin{equation}\label{ETJ}
 T_v f(x)=\partial_v f(x)+
 \sum_{\alpha\in R_+}k_\alpha\frac{\langle\alpha,v\rangle}{\langle\alpha,x\rangle}
       (f(x)-f(\sigma_\alpha x)),\qquad T_j=T_{e_j}.
\end{equation}
For $\xi\in\R^d$, the Dunkl kernel $E_k(\cdot,\xi)$ solves
$T_j E_k(\cdot,\xi)=\xi_j E_k(\cdot,\xi)$ and $E_k(0,\xi)=1$.
For $f\in L_k^1$ the Dunkl transform is defined by
\begin{equation}\label{eq:transform}
 c_k=\left(\int_{\R^d}e^{-|x|^2/2}\,d\mu_k(x)\right)^{-1},\qquad
 \Fk f(\xi)=c_k\int_{\R^d}f(x)E_k(x,-i\xi)\,d\mu_k(x).
\end{equation}
The kernel is symmetric and entire. It satisfies $E_k(gx,gy)=E_k(x,y)$ for $g\in G$ and $|E_k(ix,\xi)|\le1$ for real $x,\xi$.
The transform extends unitarily to $L_k^2$. Plancherel, inversion and the Dunkl operator identities give, for Schwartz functions~\cite{MFE1,CFD1,CFD2,CFD3,MR1},
\begin{equation}\label{eq:basic}
 \begin{gathered}
 \|\Fk f\|_{2,k}=\|f\|_{2,k},\qquad
 \Fk^{-1}f(x)=\Fk f(-x),\\
 \Fk(T_j f)(\xi)=i\xi_j\Fk f(\xi),\qquad
 T_v(f\circ g)=(T_{gv}f)\circ g.
 \end{gathered}
\end{equation}
The operators $T_j$ commute and are skew-adjoint on their natural graph
domains. The transform is a topological automorphism of the Schwartz space.
The kernel's positive integral representation~\cite{MR1} also implies
$|\partial_x^I E_k(ix,\xi)|\le |\xi|^{|I|}$ for any ordered derivative $I$.

For smooth compactly supported $f$, the Dunkl integration-by-parts formula and the expression for the Dunkl Laplacian~\cite{MR1} give
\begin{equation}\label{eq:energy}
 \sum_{j=1}^d\|T_j f\|_{2,k}^2
 =\|\nabla f\|_{2,k}^2+
 \frac12\sum_{\alpha\in R_+}k_\alpha|\alpha|^2
 \int_{\R^d}\frac{|f(x)-f(\sigma_\alpha x)|^2}{|\langle\alpha,x\rangle|^2}\,d\mu_k(x).
\end{equation}
This identity extends by closure of the nonnegative quadratic form, and also applies
to a root subsystem in its linear span with its own weight.
For $x\in\R^d$ and $r>0$, the ball-volume estimate of
\cite[Section 2]{JDA2} and \cite[(2.4)]{DHtranslations} reads
\begin{equation}\label{EWBOUND}
 \mu_k(B(x,r))\asymp r^d
 \prod_{\alpha\in R_+}(|\langle\alpha,x\rangle|+r)^{2k_\alpha}.
\end{equation}

\begin{proposition}\label{TEV}
If $E$ is $(\gamma,L)$-thick with respect to either measure, it is
$(\widetilde\gamma,L)$-thick with respect to the other, where
$\widetilde\gamma>0$ depends only on $\gamma,d$.
\end{proposition}
\begin{proof}
For another proof, see Zhao, Wei and Duan~\cite[proof of Theorem 1.1]{S4}. By homogeneity it suffices
to consider a unit cube $Q$. For $\alpha\in R_+$,
Fubini's theorem along a coordinate with
$|\alpha_j|\ge |\alpha|/\sqrt d$ gives
\[
 |\{x\in Q:|\langle\alpha,x\rangle|<\delta|\alpha|\}|
 \le2\sqrt d\,\delta.
\]
Choose $\delta>0$ so small that the strips occupy at most $\gamma/2$ of $Q$, and denote the remainder by $Q^\delta$.
For $x\in Q$, $y\in Q^\delta$,
$|\langle\alpha,x\rangle|\le|\langle\alpha,y\rangle|+|\alpha|\sqrt d$.
Consequently
\[
 \sup_Qw_k\le K_\delta\inf_{Q^\delta}w_k,
 \qquad K_\delta=(1+\sqrt d/\delta)^{2\gamma_k}.
\]
If $|E\cap Q|\ge\gamma$, then
$\mu_k(E\cap Q)\ge(\gamma/2)\inf_{Q^\delta}w_k
\ge\gamma\mu_k(Q)/(2K_\delta)$.
For the reverse implication choose a fixed $\delta$ making
$|Q^\delta|\ge1/2$. Then
$\mu_k(Q)\ge(2K_\delta)^{-1}\sup_Qw_k$, and weighted thickness implies
$|E\cap Q|\ge\gamma/(2K_\delta)$.
The case $R=\varnothing$ is immediate.
\end{proof}

Following Dziuba\'nski and Hejna~\cite[Section 2]{JDA2,DHtranslations},
we define the generalized translation on $L_k^2$ by
\begin{equation}\label{eq:translation}
 \Fk(\tau_x q)(\xi)=E_k(ix,\xi)\Fk q(\xi).
\end{equation}
Thus $\tau_x q(y)=q(y+x)$ when $k=0$. The translation denoted by
$\tau_x$ in~\cite[Section 3]{ST1} is $\tau_{-x}$ in this convention.
For $q\in\mathcal S(\R^d)$ and $f\in L_k^2$, we use the normalization
\begin{equation}\label{eq:convolution}
 (f\star_k q)(x)=c_k\int_{\R^d} f(y)\,\tau_xq(-y)\,d\mu_k(y)
 =\Fk^{-1}\!\big((\Fk f)(\Fk q)\big)(x).
\end{equation}
By inversion the corresponding product formula is
\begin{equation}\label{eq:duality}
 \Fk(fq)=(\Fk f)\star_k(\Fk q).
\end{equation}
Write $GA=\bigcup_{g\in G}gA$ and $d_G(x,y)=\min_{g\in G}|x-gy|$.
For $q\in L_k^2$ supported in $\overline B_r$, the support theorem
\cite[Theorem 1.7]{JDA2} gives
\begin{equation}\label{eq:translation-support}
 \supp\big(y\mapsto\tau_x q(-y)\big)
 \subset\bigcup_{g\in G}\overline B(gx,r).
\end{equation}
When $k_\alpha>0$ for every root, the pointwise estimate of
\cite[Theorem 4.1 and Remark 4.2]{DHtranslations} gives, for
$q\in\mathcal S(\R^d)$, $q_t(z)=t^{-d_k}q(z/t)$, $M>0$ and $t>0$,
\begin{equation}\label{eq:translated-decay}
 |\tau_x q_t(-y)|\le
 \frac{C_{q,M}}{\mu_k(B(x,t))}
 \left(1+\frac{|x-y|}{t}\right)^{-1}
 \left(1+\frac{d_G(x,y)}{t}\right)^{-M}.
\end{equation}
For radial $q(z)=\widetilde q(|z|)$ the positive representation
of~\cite{MR3,ST1} takes the form
\begin{equation}\label{eq:radial-translation}
 \tau_xq(-y)=\int_{\conv(Gx)}
 \widetilde q\!\left(\sqrt{|x|^2+|y|^2-2\langle z,y\rangle}\right)
 \,d\nu_x(z).
\end{equation}
Here $\nu_x$ is a probability measure, and the square root is at least
$d_G(x,y)$. For $k\ge0$, translation of radial functions preserves
nonnegativity and the integral and is contractive on $L_k^1$~\cite{MR3,ST1}.

\begin{lemma}\label{lem:localization}
Let $P_qf=f\star_kq$, with $q\in C_c^\infty(B_r)$.
 
If $E,S$ are measurable and $GE+B_r\subset S$, then
\begin{equation}\label{eq:physical-localization}
 \mathbf1_E P_qf=\mathbf1_E P_q(\mathbf1_Sf).
\end{equation}
If $\Psi\in\mathcal S$, $\supp\Fk\Psi\subset\overline B_r$, and
$\supp\Fk f\subset Y$, then
\begin{equation}\label{eq:spectral-localization}
 \supp\Fk(\Psi f)\subset\overline{GY+B_r}.
\end{equation}
\end{lemma}
\begin{proof}
Choose $r'<r$ with $\supp q\subset\overline B_{r'}$.
By \eqref{eq:translation-support}, for $x\in E$ the integral in
\eqref{eq:convolution} only uses $y\in Gx+\overline B_{r'}\subset S$.
For the second assertion use \eqref{eq:duality} with kernel $\Fk\Psi$.
If $\xi\notin\overline{GY+B_r}$, the kernel vanishes on $Y$ throughout
a neighborhood of $\xi$, which gives the claimed support inclusion.
\end{proof}

The spherical Paley--Wiener theorem of de Jeu~\cite[Theorem 4.10]{MFE2}
has the following form.
\begin{theorem}\label{PWSDK}
For $a>0$, $\Fk$ maps the smooth functions on $\R^d$ supported in
$\overline B_a$ bijectively onto the entire functions $f$ on $\C^d$
such that, for every integer $N\ge0$,
\begin{equation}\label{eq:PW}
 |f(z)|\le C_N(1+|z|)^{-N}e^{a|\operatorname{Im}z|}\qquad(z\in\C^d).
\end{equation}
\end{theorem}
At $k=0$ one recovers the classical Paley--Wiener theorem
\cite[Theorem 7.3.1]{LH1} in our Fourier normalization.

\subsection{Porosity and the Beurling--Malliavin multiplier theorem}
\label{subsec:classical-BM}
In the notation of Definition~\ref{def:porosity}, isometries preserve
both kinds of porosity, and dilation by $s>0$ multiplies both scale
endpoints by $s$. For line porosity, the dilation and thickening rules are given in~\cite[Lemma~3.11(a),(b)]{AC}. The same radius calculation applies to balls and to our closed scale intervals.
If $A$ is $\nu$-porous on balls or on lines from scales $\alpha_0$ to $\alpha_1$,
then for $\varepsilon>0$ its thickening $A+\overline B_\varepsilon$ is $(\nu/2)$-porous in the same sense at each scale
\begin{equation}\label{eq:porosity-thickening}
 \max\{\alpha_0,2\varepsilon/\nu\}\le r\le\alpha_1 .
\end{equation}
Indeed, at such a scale a hole of radius $\nu r$ retains its open
concentric ball of radius $\nu r-\varepsilon\ge\nu r/2$ after thickening.

For the finite unions arising from the reflection group, we record
the loss in the smallest admissible scale.
\begin{lemma}\label{lem:porosity-union}
Let $0<\nu\le1/4$ and suppose $A_1,\ldots,A_m$ are all $\nu$-porous
on balls, or all $\nu$-porous on lines, from scales $\alpha_0$ to
$\alpha_1$. Put $q=\nu/4$ and $\theta=q^m$. Their union is
$\theta$-porous in the same sense at each scale
$\alpha_0/\theta\le r\le\alpha_1$.
\end{lemma}
\begin{proof}
Start with a ball or segment of scale $r$. At step $j$ test $A_j$ on
the concentric half of the current object of scale $s_j=q^{j-1}r$.
Since $s_j/2\ge\alpha_0$, this gives a hole of radius $\nu s_j/2$.
Center the next object, of scale $q s_j$, at the hole center, on the
same line for segments. The objects are nested, and the next hole lies
inside the current hole since its center moves at most $q s_j/4$ and
its radius is $\nu q s_j/2$. The last hole has radius $2\theta r$,
its center belongs to the original object, and it misses every $A_j$.
\end{proof}

For line sections we use Cohen's estimate~\cite[Corollary~3.13]{AC}
directly. If $A$ is $\nu$-porous on lines from scales $\alpha_0$ to
$\alpha_1$ and $\tau$ is a segment of length $r\in(\alpha_0,\alpha_1)$, then
\begin{equation}\label{eq:line-length}
 |A\cap\tau|_1\le Cr(\alpha_0/r)^\zeta,\qquad C=C(\nu),\quad
 0<\zeta=\zeta(\nu)\le1 .
\end{equation}
Here $|\cdot|_1$ denotes one-dimensional Lebesgue outer measure along
the line: one may apply the estimate to $\overline A$, which has the same
open holes. At $r=\alpha_0$ use $|A\cap\tau|_1\le r$; at
$r=\alpha_1>\alpha_0$ take increasing concentric subsegments and use
$\overline A$ in \eqref{eq:line-length}.

For the higher-dimensional Beurling--Malliavin multiplier theorem we use
the formulation in Cohen's thesis~\cite[Theorem 2.5]{AC} (see also~\cite{Cohen2025}).
Given $\omega\in C^3(\R^d)$ with $\omega\le0$, put
\[
 G_\omega(x)=\int_{1/2}^2|\omega(sx)|\,ds,
 \qquad G^*(r)=\sup_{|x|=r}G_\omega(x).
\]
\begin{theorem}\label{thm:cohen-classical}
Suppose that $C_{\rm reg},C_{\rm gr}>0$ satisfy
\begin{equation}\label{eq:Cohen-hyp}
 |\partial^\beta\omega(x)|\le C_{\rm reg}\langle x\rangle^{1-|\beta|}
 \quad(|\beta|\le3),\qquad
 \int_0^\infty\frac{G^*(r)}{1+r^2}\,dr\le C_{\rm gr}.
\end{equation}
With $\sigma=C_d\max\{C_{\rm reg},C_{\rm gr}\}$ there exists
$f\in L^2(\R^d)$ such that
\begin{equation}\label{eq:Cohen-out}
 \begin{gathered}
 \supp\widehat f\subset\overline B_\sigma,\qquad
 |f(0)|\ge b(d,C_{\rm reg},C_{\rm gr})>0,\\
 |f(x)|\le e^{\omega(x)},\qquad |f(x)|\le e^{-c_d\sigma\langle x\rangle^{1/2}}
 \qquad(x\in\R^d).
 \end{gathered}
\end{equation}
\end{theorem}
The dimensional constants include the conversion from Cohen's Fourier
convention to ours. Applying the theorem to $t\omega$ multiplies the
bandwidth by $t$.

\subsection{Unique continuation weights and spectral moments}
\label{subsec:weights}
Let $W$ be a unique continuation weight as in
Definition~\ref{def:UCweight}. We use the associated sequence of spectral
moments following Jaye and Mitkovski~\cite[Section~2]{BJ}.
The model weight used below is Cohen's example with a positive constant
in the exponent:
\begin{equation}\label{eq:model-weight}
 W_{\alpha,c}(t)=\exp\!\left(\frac{ct}{(\log(2+t))^\alpha}\right),
 \qquad 0<\alpha<1,\quad c>0 .
\end{equation}
Convexity of $s\mapsto\log W_{\alpha,c}(e^s)$ is noted in~\cite[Section 3.1]{AC}. For this weight, condition~(iii) of Definition~\ref{def:UCweight} follows from the divergence of $\int_2^\infty dt/(t(\log t)^\alpha)$.

Set $D_{\eta,j}=T_j-i\eta_j$. For an ordered index
$I=(j_1,\ldots,j_m)$ write $D_{\eta,I}=D_{\eta,j_1}\cdots D_{\eta,j_m}$;
the empty index denotes the identity. Direct estimates using unshifted
Dunkl derivatives introduce dependence on $\eta$. The shifted operators
$D_{\eta,j}$ avoid this dependence in the spectral moment bounds,
as the following lemma shows.

\begin{lemma}\label{L92}
For every unique continuation weight $W$ there exists a positive sequence $(M_m)_{m\ge 0}$ such that
\begin{equation}\label{E92}
M_0 = 1,\qquad M_m \le M_{m+1},\qquad M_m^2 \le M_{m-1}M_{m+1}\ (m\ge 1),
\qquad \sum_{m=0}^{\infty}\frac{M_m}{M_{m+1}}=\infty. 
\end{equation}
If $\|W(|\cdot-\eta|)\Fk f\|_{2,k}\le C_W\|f\|_{2,k}$, then
\begin{equation}\label{E93}
\sum_{|I|=m}\|D_{\eta,I}f\|_{2,k}^2
=\bigl\||\xi-\eta|^m\Fk f\bigr\|_{2,k}^2
\le C_W^2 M_m^2\|f\|_{2,k}^2. 
\end{equation}
\end{lemma}

\begin{proof}
Suppose first that $W$ is finite. Normalize it by
\[
 W_0(t)=1\quad(0\le t\le1),\qquad
 W_0(t)=W(t)/W(1)\quad(t>1),\qquad
 M_m=\sup_{t\ge1}\frac{t^m}{W_0(t)}.
\]
Then $M_0=1$, the sequence is nondecreasing and logarithmically convex,
and $t^m\le M_mW_0(t)\le M_mW(t)$ for $t\ge0$.
The weight--sequence comparison in~\cite[Proposition 2.2]{BJ}, applied
to $W_0$, gives \eqref{E92}, since normalization changes the logarithmic
integral only by a finite quantity.
If $W(T)=\infty$ for some finite $T$, the spectral condition forces
support in $\{|\xi-\eta|\le T\}$; take $M_m=\max(1,T)^m$.
Finally, Plancherel and
\[
 \sum_{j_1,\ldots,j_m=1}^d
       \prod_{q=1}^m|\xi_{j_q}-\eta_{j_q}|^2=|\xi-\eta|^{2m}
\]
prove \eqref{E93}.
\end{proof}

\section{Center-uniform observability from Dunkl spectral decay}
\label{sec:logvinenko-sereda-main}

In this section we prove Theorem~\ref{T91}, following ideas
from~\cite{S4,BJ}. We first show that the mass in a sufficiently small
neighbourhood of the reflection hyperplanes is small, and then apply
the good-and-bad cube method to the remaining region.
 
In the first subsection, we uniformly rule out mass in a neighbourhood of the reflection hyperplanes; in the second subsection, we carry out quasi-analytic observation on the remaining orbit core; and in the third subsection, we complete the proof of the theorem. Section~\ref{subsec:logvinenko-sereda} gives the band-limited consequence, and Section~\ref{subsec:necessity} proves necessity of the logarithmic integral condition.

\subsection{Uniform non-concentration near reflection hyperplanes}
\label{sec:walls}
We give the stratified argument in full because uniformity in the frequency center is essential.

Let $H_\alpha = \{x : \langle \alpha, x\rangle = 0\}$, and set
\[
S_\delta = \bigcup_{\alpha\in R_+} \{x : \operatorname{dist}(x,H_\alpha) < \delta\},\qquad
\|D_\eta f\|_{2,k}^2 = \sum_{j=1}^d \|D_{\eta,j} f\|_{2,k}^2.
\]

In this subsection we prove a non-concentration result depending only on the first-order shifted spectral bound. It treats simultaneously intersections of hyperplanes and regions extending to infinity along the hyperplanes. 

\begin{theorem}\label{T71}
Fix $0\le B<\infty$. For any $\varepsilon>0$, there exists $\delta>0$ depending only on $\varepsilon,B,d$ such that
\begin{equation}\label{E71}
\|f\|_{2,k}=1,\quad \|D_\eta f\|_{2,k}\le B
\implies
\int_{S_\delta} |f|^2 \,d\mu_k \le \varepsilon \qquad (\eta\in\mathbb{R}^d).
\end{equation}
\end{theorem}

\begin{remark}
Lemma~5.1 of~\cite{S4} concerns spectral support contained in
$B(0,1)$. A direct application of Schur's test to translated spectral
support introduces dependence on $\eta$ through the decay estimate
in~\cite{DHtranslations}, so this argument does not yield a uniform bound.
\end{remark}

\subsubsection{Commutators and reflection overlap}
This subsubsection provides two main tools. The first is the commutator estimate \eqref{E73}, which states that after multiplication of \(f\) by a smooth cutoff function \(\chi\), the additional error incurred in the Dunkl derivatives can be controlled by \(\nabla\chi\). The second is the reflection overlap estimate \eqref{E76}: if the frequency center is large in the normal direction of some reflection hyperplane, then the local inner product between \(f\) and its reflected function can be controlled by a small quantity.

For a bounded smooth function $\chi$ with bounded gradient, expanding \eqref{ETJ} directly yields
\begin{equation}\label{E72}
  [T_v,\chi]f = (\partial_v \chi)f + \sum_{\alpha\in R_+} k_\alpha \langle \alpha, v\rangle \frac{\chi(x)-\chi(\sigma_\alpha x)}{\langle \alpha, x\rangle} f(\sigma_\alpha x). 
\end{equation}

Since $|x-\sigma_\alpha x| = 2|\langle\alpha,x\rangle|/|\alpha|$, and reflections are unitary on $L_k^2$,
\begin{equation}\label{E73}
 \|[T_v,\chi]\|_{L_k^2\to L_k^2} \le C|v|\,\|\nabla\chi\|_\infty.    
\end{equation}

Thus such multiplications preserve the first-order graph space; here $\chi$ is not required to be group-invariant or compactly supported.

Let $b=(b_j)$ be a real bounded smooth vector field with bounded first derivatives. Then
\begin{equation}\label{E74}
   \sum_j [T_j,b_j] = \operatorname{div} b + \sum_{\alpha\in R_+} k_\alpha \frac{\langle\alpha,b(x)-b(\sigma_\alpha x)\rangle}{\langle\alpha,x\rangle}\sigma_\alpha,  
\end{equation}

where $\sigma_\alpha$ on the right-hand side denotes the reflection operator acting on functions. By the skew-adjointness of $D_{\eta,j}$
\begin{equation}\label{E75}
  \left|\left\langle \sum_j [T_j,b_j]f,\,f\right\rangle_{2,k}\right|
\le 2\big\||b|f\big\|_{2,k}\,\|D_\eta f\|_{2,k}.   
\end{equation}

The multiplier $i\eta_j$ commutes with $b_j$, so there are no large-frequency terms here.

Now fix $\alpha\in R$, set $v=\alpha/|\alpha|$ and $\lambda=\langle\eta,v\rangle$. Reflection covariance gives
\[
(T_v-i\lambda)\sigma_\alpha = -\sigma_\alpha(T_v-i\lambda)-2i\lambda\sigma_\alpha.
\]

To derive the reflection-overlap estimate, write
$D_v=T_v-i\lambda=\sum_{j=1}^d v_jD_{\eta,j}$.
Since $|v|=1$, the normalization hypotheses imply
\[
\|f\|_{2,k}=1,
\qquad
\|D_v f\|_{2,k}\le \|D_\eta f\|_{2,k}\le B.
\]
The scalar multiplier $-i\lambda$ commutes with $\chi$, so
$[D_v,\chi]=[T_v,\chi]$. Combining this with the preceding
reflection identity gives
\[
\begin{aligned}
D_v(\chi\sigma_\alpha f)
&=\chi D_v\sigma_\alpha f+[T_v,\chi]\sigma_\alpha f\\
&=-\chi\sigma_\alpha D_v f
  -2i\lambda\chi\sigma_\alpha f
  +[T_v,\chi]\sigma_\alpha f.
\end{aligned}
\]
Taking the inner product with $f$, with the convention that
the inner product is linear in its first argument, and using
$D_v^*=-D_v$, we obtain
\[
\begin{aligned}
2i\lambda\langle\chi\sigma_\alpha f,f\rangle_{2,k}
={}&\langle[T_v,\chi]\sigma_\alpha f,f\rangle_{2,k}\\
&+\langle\chi\sigma_\alpha f,D_v f\rangle_{2,k}
-\langle\chi\sigma_\alpha D_v f,f\rangle_{2,k}.
\end{aligned}
\]
The measure $\mu_k$ is reflection-invariant, so $\sigma_\alpha$
is unitary on $L_k^2$. Hence the Cauchy--Schwarz inequality yields
\[
\begin{aligned}
&|\langle\chi\sigma_\alpha f,D_v f\rangle_{2,k}|
 +|\langle\chi\sigma_\alpha D_v f,f\rangle_{2,k}|\\
&\hspace{2em}\le
2\|\chi\|_\infty\|f\|_{2,k}\|D_v f\|_{2,k}
\le 2\|\chi\|_\infty B.
\end{aligned}
\]
Moreover, the commutator estimate \eqref{E73}, which follows
from \eqref{E72}, gives
\[
|\langle[T_v,\chi]\sigma_\alpha f,f\rangle_{2,k}|
\le
\|[T_v,\chi]\|_{L_k^2\to L_k^2}\|f\|_{2,k}^2
\le C\|\nabla\chi\|_\infty.
\]
No reflection invariance of $\chi$ is required here.
Taking absolute values in the preceding identity therefore gives

\begin{equation}\label{E76}
   2|\lambda|\left|\langle \chi\sigma_\alpha f,f\rangle_{2,k}\right|
\le 2\|\chi\|_\infty B + C\|\nabla\chi\|_\infty.  
\end{equation}

\subsubsection{A fixed intersection stratum}
We begin with direct estimates for a fixed stratum and a fixed outer gap, and then complete the proof using a finite covering argument.

\begin{definition}
Let $R$ be the root system and, for each $\alpha\in R$, let
\[
H_\alpha:=\{x\in\mathbb{R}^d:\langle\alpha,x\rangle=0\}.
\]
An intersection stratum is a subspace of the form
\[
\Pi=\bigcap_{\alpha\in A}H_\alpha
   =\bigl(\operatorname{span}A\bigr)^\perp,
\qquad \varnothing\ne A\subseteq R.
\]
We denote the collection of all distinct intersection strata by
\[
\mathcal{P}
:=
\left\{
\bigcap_{\alpha\in A}H_\alpha:
\varnothing\ne A\subseteq R
\right\}.
\]
In particular, every reflection hyperplane $H_\alpha$ is an
intersection stratum.
\end{definition}

An intersection stratum is a subspace obtained by intersecting a nonempty family of reflection hyperplanes. For such a stratum $\Pi$, put
\[
 R_\Pi=R\cap\Pi^\perp,\qquad x=z+y,\qquad z=P_\Pi x,\quad y=P_{\Pi^\perp}x.
\]
Here $P_V$ denotes orthogonal projection onto $V$. The roots in $R_\Pi$ form a root subsystem spanning $\Pi^\perp$, whose reflection group fixes $\Pi$ pointwise. Roots $\beta\in R_+\setminus(R_\Pi)_+$ are called outer roots. Conditions involving all outer roots are void if there are none.

Define
\[
U_{\Pi,a}(\delta) = \{x: |y|<\delta,\ |\langle \beta,z\rangle|\ge a \text{ for all outer roots } \beta\}.
\]

The main result of this subsubsection is Lemma~\ref{L72}, which states that in the region that lies close to \(\Pi\) but stays away from the other reflection hyperplanes along \(\Pi\), the mass of the function is very small.
 
Here a positive outer-root gap \(a>0\) is introduced temporarily in order to rule out locations at which a further intersection with the other hyperplanes occurs; those locations are left to Section~3.1.3.

\begin{lemma}\label{L72}
Fix $a>0$. For any $\varepsilon>0$, there exists $\delta>0$ such that all $f,\eta$ in Theorem~\ref{T71} satisfy
\begin{equation}\label{E77}
 \int_{U_{\Pi,a}(\delta)} |f|^2 d\mu_k \le \varepsilon. 
\end{equation}
\end{lemma}

\begin{proof}
\noindent\textbf{Step 1: Localization away from the outer walls.}

We construct cutoffs that separate the stratum from the outer walls,
so that the full Dunkl operators can be compared with those of $R_\Pi$.

First, take a smooth cutoff function $\psi$ on $\Pi$ with $0\le \psi\le 1$, such that $\psi(z)=1$ whenever $|\langle\beta,z\rangle|\ge a$ for all outer roots, and $\psi(z)=0$ whenever $|\langle\beta,z\rangle|\le a/2$ for some outer root, with $\|\nabla\psi\|_\infty \le C/a$.  The explicit construction is as follows.

Let
\[
q(t)=
\begin{cases}
e^{-1/t}, & t>0,\\
0, & t\le 0,
\end{cases}
\qquad
\theta(s)=\frac{q(1-s)}{q(1-s)+q(s-\tfrac12)}.
\]
Then $\theta\in C^\infty(\mathbb R)$, $0\le\theta\le1$,
$\theta=1$ on $(-\infty,1/2]$, and $\theta=0$ on $[1,\infty)$.
Define
\[
\psi(z)=
\prod_{\beta\in R_+\setminus(R_\Pi)_+}
\left[
1-\theta\!\left(\frac{|\langle\beta,z\rangle|}{a}\right)
\right],
\qquad z\in\Pi,
\]
with the empty product interpreted as $1$.
Each factor vanishes in a neighborhood of
$\{z\in\Pi:\langle\beta,z\rangle=0\}$, so the absolute value
causes no loss of smoothness. Thus $\psi\in C^\infty(\Pi)$
and $0\le\psi\le1$. Moreover, $\psi(z)=1$ whenever
$|\langle\beta,z\rangle|\ge a$ for every outer root, whereas
$\psi(z)=0$ whenever $|\langle\beta,z\rangle|\le a/2$
for some outer root. By the product rule,
\[
\|\nabla_\Pi\psi\|_\infty
\le
\frac{\|\theta'\|_\infty}{a}
\sum_{\beta\in R_+\setminus(R_\Pi)_+}|P_\Pi\beta|
\le \frac{C_{\mathrm{cut}}}{a},
\qquad
C_{\mathrm{cut}}=\|\theta'\|_\infty\sum_{\beta\in R_+}|\beta|.
\]

    We use the cutoff $\theta$ defined above, which equals $1$ on $[0,1/2]$ and $0$ on $[1,\infty)$.

Choose $0<r_0\le 1$ such that $r_0\max_{\alpha\in R}|\alpha|\le a/4$. For $z\in\operatorname{supp}\psi$ and $|y|\le r_0$, every outer root satisfies
\begin{equation}\label{E78}
   |\langle\beta,x\rangle|\ge a/4. 
\end{equation}

If $f_{\mathrm{loc}}$ is supported in this region, then
\begin{equation}\label{E79}
    \|(T_v^R - T_v^{R_\Pi})f_{\mathrm{loc}}\|_{2,k}\le \frac{C}{a}|v|\|f_{\mathrm{loc}}\|_{2,k}. 
\end{equation}

Indeed, the outer-root terms only involve $f_{\mathrm{loc}}(x)$ and $f_{\mathrm{loc}}(\sigma_\beta x)$; on the support of the latter we still have $|\langle\beta,x\rangle|=|\langle\beta,\sigma_\beta x\rangle|\ge a/4$, and reflection unitarity yields this bound. The superscripts below serve only to distinguish operators for the full root system and the subsystem.

\medskip
\noindent\textbf{Step 2: The case $\dim U>0$.}

We use the commutator and reflection-overlap estimates to bound the
mass near the stratum when the low-frequency roots do not span the
normal space.

Fix a candidate threshold $\Lambda\ge 1$, and set
\[
U = \Pi^\perp \cap \operatorname{span}\{\alpha\in R_\Pi: |\langle\alpha,\eta\rangle|/|\alpha|\le \Lambda\}^\perp.
\]

Assume $p:=\dim U>0$. For $0<r\le r_0$, let
\[
\chi_r(x)=\psi(z)\theta(|y|/r), \qquad b_r(x)=P_U x\,\chi_r(x).
\]
Since $U\subset \Pi^\perp$, we have $\|b_r\|_\infty\le r$, and
\[
\operatorname{div} b_r = p\chi_r + \rho_r,\qquad
\rho_r = \psi(z)\theta'(|y|/r)\frac{|P_U y|^2}{r|y|}.
\]
The function $\rho_r$ is uniformly bounded, supported in $r/2\le |y|\le r$, and extended by zero at $y=0$. For $\alpha\in R_\Pi$, $\chi_r$ is invariant under $\sigma_\alpha$, so we have exactly
\begin{equation}\label{E710}
   \frac{\langle\alpha,b_r(x)-b_r(\sigma_\alpha x)\rangle}{\langle\alpha,x\rangle}
=\frac{2|P_U\alpha|^2}{|\alpha|^2}\chi_r(x).  
\end{equation}

If $|\langle\alpha,\eta\rangle|/|\alpha|\le \Lambda$, this coefficient vanishes. The remaining nonvanishing terms are controlled by \eqref{E76} together with $\|\nabla\chi_r\|_\infty\le C(r^{-1}+a^{-1})$.
The outer-root terms are bounded by $Cr/a$ via \eqref{E78} and $\|b_r\|_\infty\le r$. Substituting into \eqref{E74} and \eqref{E75}, we obtain
\begin{equation}\label{E711}
   \begin{aligned}
p\int \chi_r|f|^2 d\mu_k
&\le Cr(B+a^{-1}) + C\int_{r/2\le|y|\le r}|f|^2 d\mu_k \\
&\quad + \frac{C(B+r^{-1}+a^{-1})}{\Lambda}. 
\end{aligned} 
\end{equation}
The constants here are uniform with respect to the choice of $U$, since the operator norm of the projection is at most $1$ and the number of roots is fixed.

We now choose the parameters sequentially. First pick an integer $N_0$ such that $C/N_0\le \varepsilon/4$. Next take $0<r_1\le r_0$ satisfying $Cr_1(B+a^{-1})\le \varepsilon/4$. Set
\[
r_j = r_1 4^{-(j-1)},\qquad 1\le j\le N_0,
\]
and after fixing these radii, choose $\Lambda\ge 1$ so that
\begin{equation}\label{E712}
 \frac{C(B + r_{N_0}^{-1}+a^{-1})}{\Lambda}\le \varepsilon/4.   
\end{equation}

The annuli $\{r_j/2\le |y|\le r_j\}$ are pairwise disjoint. Since $\|f\|_{2,k}^2=1$, at least one annulus has mass at most $1/N_0$. Apply \eqref{E711} for this $j$; provided the final chosen $\delta\le r_{N_0}/2$, we have $\chi_{r_j}=1$ on $U_{\Pi,a}(\delta)$, so its mass is bounded by $3\varepsilon/4$.

\medskip
\noindent\textbf{Step 3: Normal energy bounds when $U=\{0\}$.}

We bound the normal frequency and absorb the outer-root weight to
obtain a normal gradient estimate for the localized function.

Assume $U=\{0\}$. In this case the low-frequency roots span $\Pi^\perp$, hence
\begin{equation}\label{E713}
   |P_{\Pi^\perp}\eta|\le C_R\Lambda. 
\end{equation}

Indeed, $U=\{0\}$ implies that the low-frequency roots span
$\Pi^\perp$. Writing $\ell=\dim\Pi^\perp$, choose linearly
independent roots $\alpha_1,\ldots,\alpha_\ell\in R_\Pi$ such that
\[
\operatorname{span}\{\alpha_1,\ldots,\alpha_\ell\}=\Pi^\perp,
\qquad
\frac{|\langle\alpha_j,\eta\rangle|}{|\alpha_j|}
\le\Lambda
\quad(1\le j\le\ell).
\]
Set $e_j=\alpha_j/|\alpha_j|$ and $w=P_{\Pi^\perp}\eta$.
Since $e_j\in\Pi^\perp$, we have
\[
|\langle w,e_j\rangle|
=|\langle\eta,e_j\rangle|
\le\Lambda.
\]
Let $A=(\langle e_i,e_j\rangle)_{i,j=1}^{\ell}$ be the Gram
matrix of this basis. It is positive definite. Writing
$w=\sum_{j=1}^{\ell}c_je_j$ and
$b=(\langle w,e_j\rangle)_{j=1}^{\ell}$, we obtain $b=Ac$ and
\[
|w|^2
=c^{\mathsf T}Ac
=b^{\mathsf T}A^{-1}b
\le\lambda_{\min}(A)^{-1}|b|^2
\le\frac{\ell}{\lambda_{\min}(A)}\Lambda^2.
\]
Since $R$ is finite, there are only finitely many nonempty
linearly independent subsets of $R$. The Gram matrices of their
normalized roots are all positive definite, so their smallest
eigenvalues have a common positive lower bound $c_R$.
Consequently,
\[
|P_{\Pi^\perp}\eta|
\le\sqrt{\frac{\ell}{c_R}}\,\Lambda
\le C_R\Lambda,
\qquad C_R=\sqrt{\frac{d}{c_R}},
\]
which proves \eqref{E713}. The constant depends only on the fixed
root system and is uniform in $\Pi$, $\eta$, $\Lambda$, and the
choice of low-frequency root basis.

Set $f_{\mathrm{loc}} = \psi(z)\theta(|y|/r_0)f$. Define
\[
w_\Pi(y) = \prod_{\alpha\in(R_\Pi)_+} |\langle\alpha,y\rangle|^{2k_\alpha}, \qquad
w_{\mathrm{ext}}(x) = \prod_{\beta\in R_+\setminus(R_\Pi)_+} |\langle\beta,x\rangle|^{2k_\beta},
\]
so that $w_k = w_\Pi w_{\mathrm{ext}}$. The key point is that $w_{\mathrm{ext}}$ is invariant under the reflection group of $R_\Pi$: this group preserves the inner roots and outer roots separately, permits relabeling, and leaves the multiplicities unchanged.

For any unit vector $v\in\Pi^\perp$, on the support of $f_{\mathrm{loc}}$,
\begin{equation}\label{E714}
   |\partial_v \log\sqrt{w_{\mathrm{ext}}}| \le C/a.  
\end{equation}

From the commutator bound, \eqref{E79} and \eqref{E713}, for an orthonormal basis $v_1,\dots,v_\ell$ of the subspace, where $\ell=\dim\Pi^\perp$, we have
\[
\left(\sum_{j=1}^\ell \|T_{v_j}^{R_\Pi} f_{\mathrm{loc}}\|_{2,k}^2\right)^{1/2}
\le C\big(B + \Lambda + r_0^{-1} + a^{-1}\big).
\]

Absorb the outer-root weight into the function by setting $f_{\mathrm{w}}=\sqrt{w_{\mathrm{ext}}}\,f_{\mathrm{loc}}$, and let $d\mu_\Pi = w_\Pi(y)\,dy\,dz$. By subgroup invariance of the outer weight, the product rule reads
\[
T_v^{R_\Pi} f_{\mathrm{w}} = \sqrt{w_{\mathrm{ext}}}\, T_v^{R_\Pi} f_{\mathrm{loc}} + (\partial_v\sqrt{w_{\mathrm{ext}}})\,f_{\mathrm{loc}}.
\]
Combining with \eqref{E714} yields
\begin{equation}\label{E715}
    \left(\sum_{j=1}^\ell \|T_{v_j}^{R_\Pi} f_{\mathrm{w}}\|_{L^2(\mu_\Pi)}^2\right)^{1/2}
\le K := C\big(B+\Lambda+r_0^{-1}+a^{-1}\big). 
\end{equation}

The subsystem energy identity \eqref{eq:energy} gives
\begin{equation}\label{E716}
    \int |\nabla_{\Pi^\perp} f_{\mathrm{w}}|^2 d\mu_\Pi \le K^2. 
\end{equation}

We justify this application of the energy identity as follows.
As in the approximation argument at the end of the proof of
Theorem~\ref{T71} below, we first work with smooth $f$. By \eqref{E78},
\[
\operatorname{supp}f_{\mathrm{loc}}
\subset
\bigcap_{\beta\in R_+\setminus(R_\Pi)_+}
\{x:|\langle\beta,x\rangle|\ge a/4\}.
\]
Thus $\sqrt{w_{\mathrm{ext}}}$ is smooth on a neighborhood of
$\operatorname{supp}f_{\mathrm{loc}}$, even when the multiplicities are noninteger.
Moreover, $f_{\mathrm{loc}}$ vanishes in a neighborhood of every outer reflection
hyperplane. Consequently, $f_{\mathrm{w}}=\sqrt{w_{\mathrm{ext}}}\,f_{\mathrm{loc}}$ is smooth
after extension by zero across these neighborhoods.

Although $f_{\mathrm{w}}$ vanishes for $|y|\ge r_0$, it need not have compact
support in $z$. Choose $\zeta\in C_c^\infty(\Pi)$ with
$0\le\zeta\le1$ and $\zeta=1$ on $\{z\in\Pi:|z|\le1\}$, and set
\[
\zeta_n(z)=\zeta(z/n),
\qquad
f_{\mathrm{w},n}(x)=\zeta_n(P_\Pi x)f_{\mathrm{w}}(x).
\]
Then $f_{\mathrm{w},n}\in C_c^\infty(\mathbb R^d)$ and $\zeta_n(z)\to1$
pointwise. For every $v\in\Pi^\perp$ and $\alpha\in R_\Pi$,
\[
\partial_v\bigl(\zeta_n(P_\Pi x)\bigr)=0,
\qquad
P_\Pi\sigma_\alpha x=P_\Pi x.
\]
Hence both the differential and reflection terms give the exact
identities
\[
T_v^{R_\Pi}f_{\mathrm{w},n}=\zeta_nT_v^{R_\Pi}f_{\mathrm{w}},
\qquad
\nabla_{\Pi^\perp}f_{\mathrm{w},n}=\zeta_n\nabla_{\Pi^\perp}f_{\mathrm{w}}.
\]
Apply the subsystem energy identity \eqref{eq:energy} in the normal
variable $y$ and integrate in $z$. Since its reflection contribution
is nonnegative, \eqref{E715} yields
\[
\begin{aligned}
\int \zeta_n^2|\nabla_{\Pi^\perp}f_{\mathrm{w}}|^2\,d\mu_\Pi
&=\int|\nabla_{\Pi^\perp}f_{\mathrm{w},n}|^2\,d\mu_\Pi\\
&\le
\sum_{j=1}^{\ell}
\|T_{v_j}^{R_\Pi}f_{\mathrm{w},n}\|_{L^2(\mu_\Pi)}^2\\
&=
\sum_{j=1}^{\ell}
\|\zeta_nT_{v_j}^{R_\Pi}f_{\mathrm{w}}\|_{L^2(\mu_\Pi)}^2
\le K^2.
\end{aligned}
\]
Letting $n\to\infty$ and applying Fatou's lemma proves
\eqref{E716}, without assuming in advance that
$\nabla_{\Pi^\perp}f_{\mathrm{w}}\in L^2(\mu_\Pi)$.
If $\Pi=\{0\}$, the normal cutoff already makes $f_{\mathrm{w}}$ compactly
supported, so no tangential cutoff is needed.

\medskip
\noindent\textbf{Step 4: Radial non-concentration when $U=\{0\}$.}

We convert the normal gradient bound into a mass estimate on
$|y|<\delta$ that tends to zero with $\delta$.

Let
\[
D=\ell+2\sum_{\alpha\in(R_\Pi)_+}k_\alpha>0.
\]
Write $y=r\omega$ in the normal space. Fix $z,\omega$, and suppose that $f_{\mathrm{rad}}(r)=f_{\mathrm{w}}(z+r\omega)$ vanishes for $r\ge r_0$. From $f_{\mathrm{rad}}(r)=-\int_r^{2r_0}f_{\mathrm{rad}}'(s)\,ds$, the Cauchy--Schwarz inequality, and interchange of integrals,
\begin{equation}\label{E217}
    \int_0^\delta |f_{\mathrm{rad}}(r)|^2 r^{D-1}\,dr
\le I_D(\delta,2r_0)\int_0^{2r_0}|f_{\mathrm{rad}}'(s)|^2 s^{D-1}\,ds,
\end{equation}

where
\[
I_D(\delta,A)=\frac{\delta^2}{2D}+\frac{\delta^D}{D}\int_\delta^A s^{1-D}\,ds\to 0\quad(\delta\downarrow 0).
\]
Integrate over the angular and tangential variables and use \eqref{E716} to obtain
\[
\int_{|y|<\delta}|f_{\mathrm{w}}|^2 d\mu_\Pi\le K^2 I_D(\delta,2r_0).
\]
After fixing $\Lambda$ and $K$, choose
\[
0<\delta\le\min(r_{N_0}/2,r_0/2),\qquad K^2 I_D(\delta,2r_0)\le\varepsilon.
\]
On $U_{\Pi,a}(\delta)$ we have $f_{\mathrm{loc}}=f$, so this gives the desired mass bound.

We choose $\delta$ to satisfy the restrictions in both cases,
uniformly in $f$ and $\eta$. The order of parameter selection is
\[
\varepsilon,a,B \longrightarrow N_0 \longrightarrow r_1,r_{N_0} \longrightarrow \Lambda \longrightarrow K \longrightarrow \delta,
\]
and this completes the proof of the lemma. 
\end{proof}

\subsubsection{Finite covering by intersection strata}

We have not yet controlled the points that are close to \(\Pi\) and at the same time close to some outer reflection hyperplane. In this subsubsection, we resolve this problem using the following geometric fact: a point that lies close to \(\Pi\) but does not satisfy the outer-root gap condition must necessarily be close to an intersection subspace of lower dimension.

We now finish the proof of Theorem~\ref{T71}.

\begin{proof}
Let $\mathcal{P}$ be the finite collection of intersection strata defined above, and set $q=|\mathcal{P}|$. We process them in increasing order of dimension; every proper substratum is handled before the current stratum. For each intersection stratum, we take $\varepsilon/q$ as the mass target in Lemma \ref{L72}.

The minimal intersection stratum $\Pi_0=\bigcap_{\alpha\in R_+}H_\alpha$ has no outer roots, so we may take $a_{\Pi_0}=1$ and choose $\delta_{\Pi_0}>0$ by the local lemma. When treating the remaining $\Pi$, for each outer root $\beta$,
\[
\Pi_\beta=\Pi\cap H_\beta
\]
is an already processed proper substratum with $|P_\Pi\beta|>0$. First choose
\[
0<a_\Pi\le \frac18\min_{\beta\notin R_\Pi}|P_\Pi\beta|\delta_{\Pi_\beta},
\]
then select the admissible thickness under this gap via Lemma \ref{L72}, and further shrink it to $\delta_\Pi$ such that
\begin{equation}\label{E718}
    \delta_\Pi\le \frac18\min_{\beta\notin R_\Pi}\delta_{\Pi_\beta},\qquad
\int_{U_{\Pi,a_\Pi}(\delta_\Pi)}|f|^2\,d\mu_k\le \varepsilon/q.
\end{equation}
Here the minimum is taken only over positive outer roots.

If $\operatorname{dist}(x,\Pi)<\delta_\Pi$ and $x\notin U_{\Pi,a_\Pi}(\delta_\Pi)$, then there exists an outer root satisfying $|\langle\beta,P_\Pi x\rangle|<a_\Pi$. Taking the orthogonal projection onto $\Pi$, we obtain
\[
\operatorname{dist}(P_\Pi x,\Pi_\beta)=\frac{|\langle\beta,P_\Pi x\rangle|}{|P_\Pi\beta|}.
\]
Therefore
\begin{equation}\label{E719}
    \operatorname{dist}(x,\Pi_\beta)\le|x-P_\Pi x|+\operatorname{dist}(P_\Pi x,\Pi_\beta)
<\delta_\Pi+\frac{a_\Pi}{|P_\Pi\beta|}\le \frac14\delta_{\Pi_\beta}. 
\end{equation}

In other words, points violating the outer-wall gap condition lie within tubular regions of already processed deeper intersection strata. The dimension strictly decreases, so after finitely many steps one must enter some $U_{\Pi',a_{\Pi'}}(\delta_{\Pi'})$. Hence
\[
\{x:\operatorname{dist}(x,\Pi)<\delta_\Pi\}\subset \bigcup_{\Pi'\in\mathcal{P},\,\Pi'\subseteq\Pi}U_{\Pi',a_{\Pi'}}(\delta_{\Pi'}).
\]
Finally set $\delta=\min_{\alpha\in R_+}\delta_{H_\alpha}>0$, which yields
\[
S_\delta\subset\bigcup_{\Pi\in\mathcal{P}}U_{\Pi,a_\Pi}(\delta_\Pi),
\qquad
\int_{S_\delta}|f|^2\,d\mu_k\le\sum_{\Pi\in\mathcal{P}}\frac{\varepsilon}{q}=\varepsilon.
\]
This proves \eqref{E71} directly.

Lastly, for $f$ in the general first-order graph space, perform a Dunkl spectral cutoff with fixed $\eta$ to approximate $f$ and all $D_{\eta,j}f$ simultaneously; the truncated functions possess smooth representatives. Apply the proven conclusion with $B+1$, approximate termwise, normalize and remove the cutoff; since multiplication by $\mathbf{1}_{S_\delta}$ is a bounded $L^2$ operator, the mass bound passes to the limit. The chosen $\delta$ depends only on $B+1$ and remains uniform in $\eta$. Theorem \ref{T71} is proved.
\end{proof}

In particular, Lemma~\ref{L92} allows us to take $B=C_W M_1$ uniformly for all spectral classes of Theorem~\ref{T91}. Therefore, given any small fixed proportion, one may choose a thin neighborhood independent of the frequency center such that the mass contained within it does not exceed this proportion.

\subsection{Proof of weighted observability}
We combine Theorem~\ref{T71} with the quasianalytic derivative hierarchy
of Lemma~\ref{L92}. Fix the open fundamental chamber
\[
 \mathcal C=\{x:\langle\alpha,x\rangle>0\text{ for every }\alpha\in R_+\}.
\]

The goal is to convert Dunkl spectral decay into a local observation
inequality away from the reflection hyperplanes, with a constant
independent of the frequency center $\eta$. This is the content of
Lemma~\ref{lem:one-component}.

The first subsubsection constructs the conjugated orbit system,
converting Dunkl derivatives into an ordinary first-order system.
The second proves observation of the full orbit vector
(Lemma~\ref{L901}). The third proves observation of a single orbit
component (Lemma~\ref{lem:one-component}), allowing observation sets
that are not $G$-invariant.

\subsubsection{The conjugated orbit system}
We construct a first-order system for the orbit values of the shifted
Dunkl derivatives. Conjugation by the square-root weight and scalar
phases converts the weighted norms into Lebesgue norms and yields
coefficients bounded uniformly in the derivative order, the spatial
translation, and the frequency center. This system provides the
starting point for the local observation estimates below.

Take a fixed bounded domain $\Omega$ and a translation vector $x_0$. Assume $x_0+\Omega\subset\mathcal C$, and the distance from this set to all reflecting hyperplanes is at least $\delta_*>0$. For $b>0$, define the array
\begin{equation}\label{E95}
    Z_{m,g,I}(z) = b^{-1/2}\sqrt{w_k(x_0+z)}\,e^{-i\langle g^{-1}\eta,z\rangle}\bigl(D_{\eta,I}f\bigr)(g(x_0+z)),
\quad z\in\Omega,\ |I|=m.
\end{equation}

The norm of the array is taken as the sum of squares over all group indices and ordered indices. Therefore for any measurable $U\subset\Omega$, we have the exact identity
\begin{equation}\label{E96}
    \|Z_m\|_{L^2(U)}^2
=\frac1b\int_{x_0+U}\sum_{g\in G,\,|I|=m}\big|D_{\eta,I}f(gx)\big|^2\,d\mu_k(x).
\end{equation}

Let $L_i$ denote the contraction mapping from the array of the next order to the array of the current order:
\[
(L_i X)_{g,I}=\sum_{j=1}^d g_{ji}X_{g,(I,j)}.
\]
Since each $g$ is orthogonal,
\begin{equation}\label{E97}
    \sum_{i=1}^d |L_i X|^2 = |X|^2,\qquad \|L_i\|\le 1. 
\end{equation}

By definition, $T_j = D_{\eta,j}+i\eta_j$, and by direct computation (compare the kernel system in \cite{MR4}),
\[
\partial_i\sqrt{w_k(x)}
=\sqrt{w_k(x)}\sum_{\alpha\in R_+}k_\alpha\frac{\alpha_i}{\langle\alpha,x\rangle},
\]

Term by term, we obtain
\begin{equation}\label{E98}
\partial_i Z_m = L_i Z_{m+1} + B_i Z_m,   
\end{equation}
where $B_i$ acts only on the group indices, with the same matrix for each fixed $I$:
\begin{equation}\label{E99}
   (B_i Z_m)_{g,I} = \sum_{\alpha\in R_+} k_\alpha \frac{\alpha_i}{\langle \alpha,x_0+z\rangle} e^{i\langle (g\sigma_\alpha)^{-1}\eta - g^{-1}\eta,z\rangle} Z_{m,g\sigma_\alpha,I}.
\end{equation}
The phase derivative eliminates $i(g^{-1}\eta)_i$, and the square-root weight cancels the identity term from the reflection difference. The reflection acts by right multiplication $g\mapsto g\sigma_\alpha$, leaving $I$ unchanged. Each term remains a unit-modulus phase multiplied by a permutation matrix, hence
\begin{equation}\label{E910}
  \|B_i\|_{L^\infty(\Omega)} \le \gamma_k/\delta_*.   
\end{equation}

\subsubsection{Observation of the whole orbit array}
We prove a local observation estimate for the whole orbit array on
measurable sets with a fixed positive lower bound on their measure.
The argument first establishes quasianalytic uniqueness using only
bounds on the coefficients, and then uses compactness to obtain a
uniform positive observation bound. Its scalar crossing argument
follows the mechanism in \cite{FNM,BJ}; the bounded coefficients are
never differentiated.

\begin{lemma}\label{L901}
Let $\Omega \subset \mathbb{R}^d$ be a fixed bounded connected Lipschitz domain, and let $K_0, K_1 \subset \Omega$ be fixed bounded measurable sets of positive measure. Fix a sequence $M_m$ satisfying \eqref{E92}, together with $H \ge 1$, $B_0 \ge 0$, and $0 < \sigma \le |K_1|$. Suppose the array $Z_m$ satisfies \eqref{E98} in the sense of distributions, where $L_i$ is as above, $B_i\in L^\infty$ acts only on group indices with the same matrix for all orders, and
\begin{equation}\label{E911}
   \|B_i\|_{L^\infty} \le B_0,\quad \|Z_0\|_{L^2(K_0)} = 1,\quad \|Z_m\|_{L^2(\Omega)} \le H^{m+1} M_m \quad (m\ge 0). 
\end{equation}

Then there exists a constant $c>0$, depending only on the fixed parameters above, such that every measurable set $F\subset K_1$ with $|F|\ge \sigma$ satisfies
\begin{equation}\label{E912}
    \|\mathbf{1}_F Z_0\|_{L^2(\Omega)} \ge c. 
\end{equation}
\end{lemma}

\begin{proof}
    
The proof consists of four steps, ending with compactness after uniqueness on sets of positive measure.

\noindent\textbf{Step 1: Regularity of the scalar norm.}

Equations \eqref{E98} and \eqref{E97} yield

\begin{equation}\label{E913}
    \left(\sum_i |\partial_i Z_m|^2\right)^{1/2} \le |Z_{m+1}| + \sqrt{d}\,B_0 |Z_m|. 
\end{equation}

Consequently, the array at each order belongs to $H^1(\Omega)$. Let $U_m = |Z_m|$. The Sobolev chain rule gives
\begin{equation}\label{E914}
  |\nabla U_m| \le U_{m+1} + \sqrt{d}\,B_0 U_m \quad \text{almost everywhere}.   
\end{equation}

Starting from $p=2$, we apply scalar Sobolev embeddings simultaneously for all orders; each step consumes at most one higher-order $L^p$ norm. After finitely many iterations, there exist a fixed integer $s=d+2$ and constants $C_1,C_2\ge 1$ such that
\begin{equation}\label{E915}
    \|U_m\|_{L^\infty(\Omega)} \le N_m := C_1 C_2^m M_{m+s}. 
\end{equation}

For instance, when $d\ge 2$, take the finite sequence of exponents $p_j = 2d/(d-j)$, $j=0,\dots,d-1$, and finally use $W^{1,2d}\hookrightarrow L^\infty$; for $d=1$, directly apply $H^1\hookrightarrow L^\infty$. All embeddings act on the scalar quantity $U_m$, so there are no constants growing with the dimension of the array. Using \eqref{E915} at orders $m,m+1$ in \eqref{E913}, each $Z_m$ and $U_m$ admits a locally Lipschitz representative. The sequence $N_m$ is nondecreasing, logarithmically convex, and
\begin{equation}\label{E916}
    \sum_{m=0}^{\infty}\frac{N_m}{N_{m+1}}
= C_2^{-1}\sum_{m=0}^{\infty}\frac{M_{m+s}}{M_{m+s+1}}
=\infty. 
\end{equation}

\noindent\textbf{Step 2: Zero sets of positive measure force vanishing at all orders.} 

Suppose $Z_0=0$ on a set $A\subset\Omega$ of positive measure. The locality property of Sobolev functions implies that if $Z_m=0$ almost everywhere on $A$, then all its weak derivatives also vanish on $A$. From \eqref{E98}, $L_i Z_{m+1}=0$ on $A$, and \eqref{E97} then yields $Z_{m+1}=0$ on $A$. By induction, all orders vanish almost everywhere on the same set $A$. Take a Lebesgue density point $x_*$ of $A$ inside $\Omega$; using continuous representatives at each order and writing the exceptional null set as a countable union, we have
\begin{equation}\label{E917}
    U_m(x_*) = 0 \quad (m\ge 0). 
\end{equation}

\noindent\textbf{Step 3: Quasi-analytic propagation along finite paths.} 

Take any $x\in\Omega$, and connect $x_*$ to $x$ by a finite polygonal path compactly contained in $\Omega$, parametrized by arc length as $\Gamma:[0,T]\to\Omega$. Restricting \eqref{E913} to this path requires some justification: mollify by convolution on a compact neighborhood of each line segment, integrate the gradient bound along the segment, then let the mollification scale tend to zero. Since $U_m$ and the right-hand side are continuous, we obtain
\[
|U_m(\Gamma(t)) - U_m(\Gamma(s))| \le \int_{s}^{t} \big(U_{m+1}(\Gamma(r)) + \sqrt{d}\,B_0 U_m(\Gamma(r))\big)\,dr.
\]

This inequality holds for $0\le s\le t\le T$ and remains valid after concatenating line segments. Thus $f_m(t):=U_m(\Gamma(t))$ is absolutely continuous, and for a fixed constant $C\ge 1$,
\begin{equation}\label{E918}
    0\le f_m \le N_m,\quad |f_m'|\le C(f_m + f_{m+1}),\quad f_m(0)=0. 
\end{equation}

We now present the one-dimensional crossing argument needed. Let
\[
a_m(t)=\frac{f_m(t)}{e^m N_m},\qquad b(t)=\sup_{m\ge 0}a_m(t),\qquad r_m=\frac{N_{m+1}}{N_m}.
\]
Since $a_m\le e^{-m}$, the tails vanish uniformly, so $b$ is continuous and its positive maximum is attained at a finite index. On intervals where $b\ge e^{-q}$, the index achieving the supremum satisfies $m\le q$, and
\[
|a_m'|\le C(a_m + e r_m a_{m+1})\le C(1+e r_q)\,b \quad (m\le q).
\]
The maximum of finitely many absolutely continuous functions is still absolutely continuous, hence
\begin{equation}\label{E919}
   |(\log b)'|\le C(1+e r_q) \quad \text{almost everywhere on the set where }b\ge e^{-q}. 
\end{equation}

If $b(T)>0$, choose $q_0$ such that $e^{-(q_0-1)}<b(T)$. For $q\ge q_0-1$, let $t_q$ be the first time at which $e^{-q}$ is attained; for $q\ge q_0$, let $s_q$ be the last time before $t_{q-1}$ at which $e^{-q}$ is attained. Then
\[
t_q \le s_q < t_{q-1},\qquad e^{-q}<b(t)<e^{-(q-1)} \quad (s_q<t<t_{q-1}).
\]
By \eqref{E919}, each crossing consumes a length of at least $1/\big[C(1+e r_q)\big]$. These intervals are pairwise disjoint in their interiors, so
\[
T\ge \sum_{q\ge q_0}(t_{q-1}-s_q)\ge \frac{1}{C(1+e)}\sum_{q\ge q_0}\frac{N_q}{N_{q+1}}=\infty.
\]

\noindent\textbf{Step 4: Attaining a strictly positive minimum over a compact class.} 

Relax the indicator function of the observation set to an $L^\infty$ function $\theta$ satisfying
\[
0\le \theta\le 1,\quad \theta=0 \text{ on } \Omega\setminus K_1,\quad \int_\Omega \theta\,dx \ge \sigma.
\]
Consider the collection of all triples $(Z_m,B_i,\theta)$ satisfying these conditions together with the hierarchy \eqref{E98} and the bounds \eqref{E911}. If this class is empty, the lemma holds trivially; we assume hereafter that it is nonempty.

By \eqref{E913} and \eqref{E911}, for each fixed order $m$, $Z_m$ is uniformly bounded in $H^1(\Omega)$. For each fixed $m$, the arrays $Z_m^{(n)}$ take values in a
fixed finite-dimensional space and are uniformly bounded in
$H^1(\Omega)$. Rellich compactness and a diagonal selection
therefore yield a common subsequence, still indexed by $n$, such that
\[
Z_m^{(n)}\longrightarrow Z_m
\quad\text{strongly in }L^2(\Omega)
\qquad\text{for every }m\ge0.
\]
Here convergence is asserted at each fixed order; no uniformity
in $m$ is required. Since $L^1(\Omega)$ is separable, the uniform
$L^\infty$ bounds also allow us, after passing to a further
subsequence, to arrange simultaneously that
\[
B_i^{(n)}\rightharpoonup^* B_i
\quad\text{in }L^\infty(\Omega),\qquad i=1,\ldots,d,
\qquad
\theta^{(n)}\rightharpoonup^*\theta
\quad\text{in }L^\infty(\Omega).
\]
The matrix convergence is understood entrywise on the group
indices. At every order, the same matrix acts on these indices
and leaves the ordered derivative index unchanged; this
structure is preserved in the limit.

To pass to the limit in the coefficient terms, fix $i,m$ and write
\[
B_i^{(n)}Z_m^{(n)}-B_iZ_m
=
B_i^{(n)}(Z_m^{(n)}-Z_m)
+(B_i^{(n)}-B_i)Z_m.
\]
The first term tends to zero strongly in $L^2(\Omega)$, since
\[
\bigl\|B_i^{(n)}(Z_m^{(n)}-Z_m)\bigr\|_{L^2(\Omega)}
\le B_0\|Z_m^{(n)}-Z_m\|_{L^2(\Omega)}
\longrightarrow0.
\]
For any fixed $L^2$ test array $\varphi$, each scalar product
of a component of $Z_m$ with a component of
$\overline{\varphi}$ belongs to $L^1(\Omega)$.
The weak-$*$ convergence of the matrix entries consequently gives
\[
\int_\Omega
\bigl\langle (B_i^{(n)}-B_i)Z_m,\varphi\bigr\rangle\,dx
\longrightarrow0.
\]
Thus $B_i^{(n)}Z_m^{(n)}\rightharpoonup B_iZ_m$ weakly in
$L^2(\Omega)$, and in particular in the sense of distributions.
Passing to the limit in the distributional formulation of
\eqref{E98}, using also the strong convergence of
$Z_{m+1}^{(n)}$, yields
\[
\partial_iZ_m=L_iZ_{m+1}+B_iZ_m
\qquad(m\ge0,\ i=1,\ldots,d).
\]

The bounds in \eqref{E911} remain valid: strong $L^2$ convergence
preserves the array norm bounds, while weak-$*$ lower
semicontinuity gives $\|B_i\|_{L^\infty}\le B_0$.
Moreover,
\[
\|Z_0^{(n)}-Z_0\|_{L^2(K_0)}
\le \|Z_0^{(n)}-Z_0\|_{L^2(\Omega)}
\longrightarrow0,
\]
so the core normalization passes to the limit:
$\|Z_0\|_{L^2(K_0)}=1$.

Finally, testing the weak-$*$ convergence of $\theta^{(n)}$
against nonnegative $L^1$ functions shows that
$0\le\theta\le1$ almost everywhere. Testing against
$\mathbf1_{\Omega\setminus K_1}$ and $1$, respectively, gives
\[
\int_{\Omega\setminus K_1}\theta\,dx=0,
\qquad
\int_\Omega\theta\,dx
=
\lim_{n\to\infty}\int_\Omega\theta^{(n)}\,dx
\ge\sigma.
\]
Since $\theta\ge0$, the first identity implies that
$\theta=0$ almost everywhere on $\Omega\setminus K_1$.
Hence the limiting arrays, coefficients, and observation weight
satisfy all the defining conditions of the admissible class.

Take a minimizing sequence from the above class. Since $|Z_0^{(n)}|^2\to |Z_0|^2$ strongly in $L^1$, the functional
\[
\int_\Omega \theta |Z_0|^2\,dx
\]
is continuous along this convergence. Hence it attains its minimum over the class. If the minimum were zero, then $Z_0=0$ on $\{\theta>0\}$, and this set has measure at least $\sigma$. The uniqueness result from the preceding three steps forces $Z_0\equiv 0$, contradicting the core normalization $\|Z_0\|_{L^2(K_0)}=1$. Therefore the minimum equals some $c^2>0$. Taking $\theta=\mathbf{1}_F$ recovers \eqref{E912}.
\end{proof}

\subsubsection{Observation of one orbit component}
We prove a uniform local observation estimate for a single orbit
component carrying a fixed positive proportion of the core mass.
Using the specific coefficients in \eqref{E99}, we identify analytic
coefficients in the limiting system and establish scalar
quasianalyticity of its components. This estimate allows us to select
a component carrying maximal mass, as in
\cite[proof of Lemma 5.7]{S4}, without imposing invariance on the
observation set.

\begin{lemma}\label{lem:one-component}
Use the fixed data of Lemma~\ref{L901}, and fix $\delta_*>0$. Suppose \eqref{E98} and
\eqref{E911} hold, with $B_i$ given by \eqref{E99} for some $x_0,\eta$
such that $x_0+\Omega\subset\mathcal C$ has distance at least
$\delta_*$ from every reflection hyperplane. If $g_0\in G$ satisfies
\[
 \|Z_{0,g_0}\|_{L^2(K_0)}^2\ge |G|^{-1},
\]
then there is $c>0$, depending only on the fixed data and $\delta_*$,
such that
\[
 \|\mathbf1_F Z_{0,g_0}\|_{L^2(\Omega)}\ge c
 \qquad(F\subset K_1,\ |F|\ge\sigma).
\]
In particular, $c$ is independent of $x_0,\eta,g_0$ and $F$.
\end{lemma}
\begin{proof}
\noindent\textbf{Step 1: Compactness and the limiting coefficients.}

We identify the weak limits of the coefficients before using ordinary
derivatives. Suppose the assertion fails. Since $G$ is finite, there
are a fixed $g_0$ and sequences of admissible arrays $Z_m^{(n)}$,
parameters $x_n,\eta_n$, and sets $F_n$ such that
\[
 \|\mathbf1_{F_n}Z_{0,g_0}^{(n)}\|_{L^2}\longrightarrow0.
\]
As in Step 4 of Lemma~\ref{L901}, \eqref{E913} and Rellich compactness
allow a diagonal subsequence with
\[
 Z_m^{(n)}\longrightarrow Z_m\quad\hbox{strongly in }L^2(\Omega)
 \quad\text{for every }m,\qquad
 \mathbf1_{F_n}\rightharpoonup^*\theta\quad\hbox{in }L^\infty(\Omega).
\]
The bounds \eqref{E911}, the normalization, and the lower bound on the
$g_0$ component pass to the limit. Also $0\le\theta\le1$,
$\theta=0$ outside $K_1$, and $\int\theta\ge\sigma$.

For each positive root, pass to a further subsequence such that
$a_{\alpha,n}=\langle\alpha,x_n\rangle$ has a finite limit
$a_\alpha$ or tends to $+\infty$. These sequences are bounded below,
since $\Omega$ is bounded and
$\langle\alpha,x_n+z\rangle\ge |\alpha|\delta_*$ on $\Omega$.
The amplitudes
\[
 A_{\alpha,i}^{(n)}(z)
 =\frac{k_\alpha\alpha_i}{a_{\alpha,n}+\langle\alpha,z\rangle}
\]
therefore converge uniformly on $\Omega$, with every derivative, to
$A_{\alpha,i}(z)=k_\alpha\alpha_i/(a_\alpha+\langle\alpha,z\rangle)$
in the finite case, and to zero otherwise.
For the finitely many pairs $(g,\alpha)$, we may simultaneously arrange
that
\[
 \lambda_{g,\alpha,n}
 =(g\sigma_\alpha)^{-1}\eta_n-g^{-1}\eta_n
\]
either converges to a finite vector or has norm tending to infinity.
In the former case
$A_{\alpha,i}^{(n)}e^{i\langle\lambda_{g,\alpha,n},z\rangle}$
converges uniformly to an analytic function. In the latter case it
converges weakly-$*$ to zero in $L^\infty(\Omega)$: test against an
$L^1$ function, replace the amplitude by its uniform limit, and apply
the Riemann--Lebesgue lemma.

Thus $B_i^{(n)}\rightharpoonup^* B_i^\infty$, where $B_i^\infty$ is
analytic and still acts only on group indices, identically at every
order. Strong $L^2$ convergence of $Z_m^{(n)}$ permits passage to the
limit in the products, giving
\[
 \partial_iZ_m=L_iZ_{m+1}+B_i^\infty Z_m.
\]
The separated denominators and the finitely many remaining frequencies
give constants $A_0,A_1\ge1$ such that
\[
 \|\partial^\beta B_i^\infty\|_{L^\infty(\Omega)}
 \le A_0A_1^{|\beta|}|\beta|!.
\]
These constants may depend on the subsequence limit. No uniform
bounds on derivatives of the original $B_i^{(n)}$ are asserted or used.

\medskip
\noindent\textbf{Step 2: Scalar quasianalyticity of the limit.}

We turn the spectral hierarchy into ordinary derivative bounds for
each component of the limiting array. Step 1 of Lemma~\ref{L901},
applied to the limiting system, gives
\[
 \|Z_m\|_{L^\infty(\Omega)}\le A H_1^m M_{m+s}
\]
for fixed $s\in\mathbb N$ and $A,H_1\ge1$. Define
\[
 \mu_n=\frac{M_n}{M_{n-1}},\qquad
 \widetilde\mu_n=\max(\mu_n,n),\qquad
 \widetilde M_0=1,\qquad
 \widetilde M_n=\prod_{r=1}^n\widetilde\mu_r.
\]
The sequences of ratios are nondecreasing, $M_n\le\widetilde M_n$,
and $\widetilde\mu_n\ge n$. Moreover,
\[
 \sum_{n\ge1}\frac1{\widetilde\mu_n}=\infty.
\]
Indeed, Cauchy's condensation test and \eqref{E92} give
$\sum_{j\ge0}2^j/\mu_{2^j}=\infty$. Truncating the terms of a
divergent positive series at $1$ preserves divergence, and
\[
 \frac{2^j}{\widetilde\mu_{2^j}}
 =\min\left(\frac{2^j}{\mu_{2^j}},1\right).
\]
A second application of the condensation test proves the assertion.

Choose $D\ge\max(2H_1,2A_1,4A_0)$. Differentiating the limiting
system inductively gives
\[
 \|\partial^\beta Z_m\|_{L^\infty(\Omega)}
 \le A H_1^m D^{|\beta|}\widetilde M_{m+|\beta|+s}.
\]
For completeness, let $n=|\beta|$. The term
$L_i\partial^\beta Z_{m+1}$ costs at most $H_1/D$ times the proposed
bound of order $n+1$. In the Leibniz sum, collect all derivatives
of $B_i^\infty$ of order $l$. The identities
$\sum_{\zeta\le\beta,\,|\zeta|=l}\binom\beta\zeta=\binom nl$
and $\widetilde\mu_r\ge r$ imply
\[
 \frac{n!}{(n-l)!}\,
 \frac{\widetilde M_{m+n-l+s}}{\widetilde M_{m+n+s+1}}
 \le\frac1{n+1}.
\]
Hence the Leibniz terms cost at most
\[
 \frac{A_0}{D(n+1)}\sum_{l=0}^n(A_1/D)^l
 \le\frac{2A_0}{D}\le\frac12,
\]
which closes the induction. This argument also proves interior
smoothness, starting from the locally Lipschitz representatives in
Step 1 of Lemma~\ref{L901} and differentiating in distributions.

In particular, each component $Z_{0,g}$ has ordinary derivative bounds
with the logarithmically convex sequence $\widetilde M_{n+s}$, and
\[
 \sum_{n\ge0}\frac{\widetilde M_{n+s}}{\widetilde M_{n+s+1}}=\infty.
\]
The Denjoy--Carleman criterion, or the crossing argument in Step 3 of
Lemma~\ref{L901} applied along line segments, therefore makes each
component quasianalytic on the connected domain $\Omega$; compare
\cite{FNM,BJ}.

\medskip
\noindent\textbf{Step 3: Contradiction from observation of one component.}

Strong $L^2$ convergence and weak-$*$ convergence of the indicators give
\[
 \int_\Omega\theta|Z_{0,g_0}|^2\,dz=0.
\]
The set $\{\theta>0\}$ has measure at least $\sigma$, and the component
vanishes almost everywhere there. Sobolev locality, applied successively
to its smooth derivatives, shows that all derivatives vanish almost
everywhere on that same set. At a density point, continuity makes all
these derivatives zero. Quasianalyticity and connectedness imply
$Z_{0,g_0}\equiv0$, contradicting
$\|Z_{0,g_0}\|_{L^2(K_0)}^2\ge |G|^{-1}$. This proves the uniform
positive lower bound.
\end{proof}

\subsection{Proof of Theorem~\ref{T91}}
In this subsection, we finish the proof of Theorem~\ref{T91}.

\begin{proof}[Proof of Theorem~\ref{T91}]

We now sum directly over every function satisfying the hypotheses of Theorem~\ref{T91}. We may assume $\|f\|_{2,k}=1$. By Lemma~\ref{L92}, $\|D_\eta f\|_{2,k}\le C_W M_1$, where $\|D_\eta f\|_{2,k}^2=\sum_j\|D_{\eta,j}f\|_{2,k}^2$. By Theorem \ref{T71}, we may choose $\delta>0$ depending only on fixed data such that
\begin{equation}\label{E81}
    \|\mathbf{1}_{\mathbb{R}^d\setminus S_{2\delta}}f\|_{2,k}^2\ge \frac34. 
\end{equation}

If $R=\emptyset$, we simply set $S_{2\delta}=\emptyset$, and all subsequent geometric requirements concerning walls may be omitted.

\subsubsection*{Fixed cores and observation boxes}

Away from the reflection hyperplanes, $\{g\mathcal{C}:g\in G\}$ forms a disjoint partition. Therefore
\begin{equation}\label{E82}
    \int_{\mathcal{C}}\sum_{g\in G}|f(gx)|^2\,d\mu_k(x)=\|f\|_{2,k}^2. 
\end{equation}

Fix $\rho>0$ such that $\sqrt{d}\rho<\delta/4$. In the axis-parallel lattice tiling of side length $\rho$, we retain all cubes $x_j+Q_\rho$ whose centers $x_j$ satisfy
\[
\frac{\langle\alpha,x_j\rangle}{|\alpha|}\ge\delta \qquad (\alpha\in R_+).
\]
These cubes have pairwise disjoint interiors, are contained in $\mathcal{C}$, and cover $\mathcal{C}\setminus S_{2\delta}$ up to null boundaries.

By Proposition~\ref{TEV}, $E$ is Lebesgue $(\gamma_0,L)$-thick for
some fixed $\gamma_0>0$. Put $\ell=\sqrt d\,L$, so every orthogonal
image of $Q_\ell$ contains the axis-parallel cube $Q_L$.
Choose $v\in\mathcal{C}$, and then take a sufficiently large fixed $t$ such that the signed distance from $tv+Q_\ell$ to every wall is greater than $1$. Slightly thickening the convex hull of $Q_\rho$ and $tv+Q_\ell$, we obtain a fixed bounded convex Lipschitz domain $\Omega$ satisfying
\begin{equation}\label{E83}
   Q_\rho\cup(tv+Q_\ell)\subset\Omega,\qquad
\frac{\langle\alpha,z\rangle}{|\alpha|}>-\delta/2 \quad (z\in\Omega,\,\alpha\in R_+).   
\end{equation}

Hence each $x_j+\Omega$ lies inside $\mathcal{C}$ and has distance at least $\delta/2$ to all walls. Since the lattice spacing is fixed and $\Omega$ is bounded, the overlap multiplicity of these translated domains is bounded by a fixed integer $N$. The same $\Omega$ is used no matter how far the spatial location $x_j$ lies.

\subsubsection*{Good cores at every derivative order}

For each retained core cube define
\[
m_j=\int_{x_j+Q_\rho}\sum_{g\in G}|f(gx)|^2\,d\mu_k(x),\qquad
b_{j,m}=\int_{x_j+\Omega}\sum_{g\in G,\,|I|=m}|D_{\eta,I}f(gx)|^2\,d\mu_k(x).
\]
From \eqref{E81}, \eqref{E82} and finite overlap,
\begin{equation}\label{E84}
   \frac34\le\sum_j m_j\le 1,\qquad
\sum_j b_{j,m}\le N C_W^2 M_m^2 \quad (m\ge 0).  
\end{equation}

A core cube is called good if $m_j>0$ and for every $m\ge 0$,
\begin{equation}\label{E85}
    b_{j,m}\le 2^{m+4} N C_W^2 M_m^2 m_j. 
\end{equation}

The remaining positive-mass cores are called bad cubes; zero-mass cores may be discarded directly. By \eqref{E84}, the total mass of cores violating the order-$m$ condition is at most $2^{-m-4}$. Taking the union over all orders, the total mass of bad cores is at most $1/8$, so
\begin{equation}\label{E86}
    \sum_{j\text{ good}}m_j\ge \frac58. 
\end{equation}

Let $P_j=x_j+tv+Q_\ell$. For every $g\in G$, the set $gP_j$
contains an axis-parallel cube of side $L$. Consequently
$|E\cap gP_j|\ge\gamma_0L^d$. Define
\begin{equation}\label{E87}
 \begin{gathered}
 F_{j,g}=\{z\in tv+Q_\ell:g(x_j+z)\in E\},\\
 |F_{j,g}|=|E\cap gP_j|\ge\gamma_0L^d=:\sigma>0.
 \end{gathered}
\end{equation}

\subsubsection*{Summation of the local estimates}

For each good core, take $x_0=x_j$, $b=m_j$ in \eqref{E95} to obtain the array on the fixed domain $\Omega$. By \eqref{E96} and \eqref{E85},
\[
\|Z_0\|_{L^2(Q_\rho)}=1,\qquad
\|Z_m\|_{L^2(\Omega)}\le 4\sqrt{N}C_W(\sqrt{2})^m M_m.
\]
Thus we may choose a fixed $H\ge\max(1,4\sqrt N C_W,\sqrt2)$
satisfying \eqref{E911}. Choose $g_j\in G$ so that
\[
 \int_{x_j+Q_\rho}|f(g_jx)|^2\,d\mu_k(x)
 =\max_{g\in G}\int_{x_j+Q_\rho}|f(gx)|^2\,d\mu_k(x)
 \ge\frac{m_j}{|G|}.
\]
Then $\|Z_{0,g_j}\|_{L^2(Q_\rho)}^2\ge |G|^{-1}$.
Apply Lemma~\ref{lem:one-component} with $K_0=Q_\rho$,
$K_1=tv+Q_\ell$, $\delta_*=\delta/2$, and $F=F_{j,g_j}$.
The bound in \eqref{E910} allows $B_0=2\gamma_k/\delta$.
The lemma and \eqref{E87} give a constant $c>0$, independent of
$j,\eta,f,E$, such that
\begin{equation}\label{E88}
 \begin{aligned}
 c^2m_j
 &\le m_j\int_{F_{j,g_j}}|Z_{0,g_j}(z)|^2\,dz\\
 &=\int_{P_j\cap g_j^{-1}E}|f(g_jx)|^2\,d\mu_k(x)
 =\int_{g_jP_j\cap E}|f(y)|^2\,d\mu_k(y).
 \end{aligned}
\end{equation}
The last equality uses invariance of the measure, not of $E$.
For each fixed $g$, the family $\{gP_j\}_j$ has overlap at most $N$.
Moreover $gP_j\subset g\mathcal C$, and the chambers for different
$g$ are disjoint. Hence the selected family $\{g_jP_j\}_j$ also has
overlap at most $N$. Together with \eqref{E86}, this yields
\[
 \begin{aligned}
 N\|\mathbf1_E f\|_{2,k}^2
 &\ge\sum_{j\text{ good}}\int_{g_jP_j\cap E}|f|^2\,d\mu_k\\
 &\ge c^2\sum_{j\text{ good}}m_j\ge\frac58c^2.
 \end{aligned}
\]
Restoring homogeneity gives
\begin{equation}\label{E89}
 \|f\|_{2,k}\le\frac1c\sqrt{\frac{8N}{5}}\,
                      \|\mathbf1_E f\|_{2,k}.
\end{equation}

This proves Theorem \ref{T91}.
\end{proof}

\subsection{A Logvinenko--Sereda theorem uniform in the spectral center}
\label{subsec:logvinenko-sereda}
The band-limited formulation in \cite[Theorem 1.3]{S4} concerns a spectral ball centered at the origin. We now
state its center-uniform counterpart for arbitrary thick observation
sets, as in Theorem~\ref{T91}. The bandwidth and the center are separate
parameters: the estimate below is uniform in the center, while its
constant may depend on the bandwidth.

\begin{theorem}
\label{thm:logvinenko-sereda}
Fix $0<\gamma\le1$, $L>0$, and $N>0$. There is a finite constant
$C=C(\gamma,LN,d)$ such that, for every
measurable set $E\subset\R^d$ that is
$(\gamma,L)$-thick for $\mu_k$, every $\eta\in\R^d$, and every
$f\in L_k^2(\R^d)$,
\begin{equation}\label{eq:logvinenko-sereda-center}
 \supp\Fk f\subset\overline{B(\eta,N)}
 \quad\Longrightarrow\quad
 \|f\|_{2,k}\le C(\gamma,LN,d)
                      \|\mathbf1_E f\|_{2,k}.
\end{equation}
The dependence on the fixed Dunkl structure is suppressed as in
Section~\ref{subsec:notation}. In particular, the constant is independent
of $\eta$, of $f$, and of the particular set $E$ within this thickness
class. 
\end{theorem}

\begin{proof}
First consider $L=1$ and take the weight $W$ in \eqref{eq:example}.
It satisfies Definition~\ref{def:UCweight}, with convexity understood
in the extended-valued sense. The spectral support condition gives
\begin{equation}\label{eq:logvinenko-sereda-spectral-weight}
 \big\|W(|\cdot-\eta|)\Fk f\big\|_{2,k}=\|\Fk f\|_{2,k}=\|f\|_{2,k},
\end{equation}
where $0\cdot\infty=0$. Theorem~\ref{T91} with $C_W=1$ proves
\eqref{eq:logvinenko-sereda-center} at this scale, uniformly in $\eta$.

For general $L$, use the unitary dilation
\[
 f_L(x)=L^{d_k/2}f(Lx),\qquad E_L=L^{-1}E.
\]
A change of variables in \eqref{eq:transform} gives
\begin{equation}\label{eq:logvinenko-sereda-scaling}
 \begin{gathered}
 \Fk f_L(\xi)=L^{-d_k/2}\Fk f(\xi/L),\qquad
 \supp\Fk f_L\subset\overline{B(L\eta,LN)},\\
 \|f_L\|_{2,k}=\|f\|_{2,k},\qquad
 \|\mathbf1_{E_L}f_L\|_{2,k}=\|\mathbf1_E f\|_{2,k}.
 \end{gathered}
\end{equation}
Moreover $E_L$ is $(\gamma,1)$-thick for $\mu_k$:
for every $x\in\R^d$,
\[
 \mu_k\big(E_L\cap(x+Q_1)\big)
 =L^{-d_k}\mu_k\big(E\cap(Lx+Q_L)\big)
 \ge\gamma\mu_k(x+Q_1).
\]
Applying the unit-scale estimate to $f_L$ gives dependence only on $\gamma,LN,d$ and the fixed Dunkl structure.
Finally, Proposition~\ref{TEV} converts Lebesgue thickness into
$\mu_k$-thickness without introducing any dependence on $L$ or $\eta$.
\end{proof}

\begin{corollary}\label{cor:logvinenko-sereda-diameter}
Under the observation hypotheses of Theorem~\ref{thm:logvinenko-sereda}, let
$K\subset\R^d$ be nonempty and measurable, with
$\operatorname{diam}K\le D$ for some $D>0$. Then
\[
 \supp\Fk f\subset K
 \quad\Longrightarrow\quad
 \|f\|_{2,k}\le C(\gamma,LD,d)
                         \|\mathbf1_E f\|_{2,k}.
\]
The estimate is uniform over the position and shape of $K$.
\end{corollary}
\begin{proof}
Choose any $\eta\in K$. Then $K\subset\overline{B(\eta,D)}$, so
Theorem~\ref{thm:logvinenko-sereda} applies with $N=D$.
\end{proof}

\subsection{Necessity of the logarithmic integral condition}
\label{subsec:necessity}
Following~\cite{BJ}, we prove necessity of the logarithmic integral
condition in Definition~\ref{def:UCweight}.

\begin{theorem}\label{thm:necessity}
Let $W:[0,\infty)\to[1,\infty]$ be nondecreasing and lower
semicontinuous, with $W(0)=1$, satisfying conditions~(i) and~(ii)
of Definition~\ref{def:UCweight}. Suppose that
\[
 \int_0^\infty\frac{\log W(t)}{1+t^2}\,dt<\infty.
\]
For every $r>0$, there exists a finite constant $C_W\ge1$ such
that, for every $\eta\in\R^d$, one can find
$f_\eta\in C_c^\infty(B_r)$ satisfying
\[
 \|f_\eta\|_{2,k}=1,\qquad
 \bigl\|W(|\cdot-\eta|)\Fk f_\eta\bigr\|_{2,k}\le C_W.
\]
Consequently, for every $0<\gamma<1$ and $L>0$, there exist a
fixed $(\gamma,L)$-thick set $E$ for $\mu_k$ and a finite
$C_W\ge1$ such that, for every $\eta\in\R^d$, there is a normalized
$f_\eta\in C_c^\infty(\R^d)$ with
\[
 \bigl\|W(|\cdot-\eta|)\Fk f_\eta\bigr\|_{2,k}\le C_W,
 \qquad \mathbf1_Ef_\eta=0.
\]
Both $E$ and $C_W$ are independent of $\eta$.
Thus condition~(iii) is necessary and sufficient for the conclusion
of Theorem~\ref{T91} to hold for every $C_W,\gamma,L$, uniformly
in the spectral center.
\end{theorem}

\begin{proof}
We adapt the Paley--Wiener construction of
Jaye and Mitkovski~\cite[Section~4]{BJ}.
The finite logarithmic integral forces $W$ to be finite everywhere.
Fix $r>0$, choose $m>d_k/2$, and set
\[
 V(t)=W(\sqrt d\,t)(1+t^2)^{m/2}.
\]
This weight has a finite logarithmic integral and satisfies the
same monotonicity, convexity, and growth assumptions.
The one-dimensional sinc-product construction in
\cite[Section~4]{BJ}, followed by normalization at zero, gives
$u\in C_c^\infty((-r/(2\sqrt d),r/(2\sqrt d)))$ such that
\[
 \widehat u(0)=1,\qquad
 |\widehat u(t)|\le C\,V(|t|)^{-1}\qquad(t\in\R).
\]
Indeed, the summable moment ratios make the support arbitrarily
small, the sinc product is nonzero at zero, and its
superpolynomial decay gives smoothness.
Set
\[
 U(x)=\prod_{j=1}^d u(x_j),\qquad
 F(\xi)=\widehat U(\xi)=\prod_{j=1}^d\widehat u(\xi_j).
\]
Since $\max_j|\xi_j|\ge|\xi|/\sqrt d$ and
$\prod_j(1+|\xi_j|^2)\ge1+|\xi|^2$, we have
\[
 F(0)=1,\qquad
 |F(\xi)|\le\frac{A}{W(|\xi|)\langle\xi\rangle^m}.
\]
Fix $0<\rho\le1$ such that $|F|\ge1/2$ on $B_\rho$.

The ordinary Paley--Wiener theorem applied to
$U\in C_c^\infty(B_{r/2})$ gives \eqref{eq:PW} for $F$
with radius $r/2$.
For every $\eta\in\R^d$ and integer $N\ge0$,
\[
 |F(z-\eta)|
 \le A_N(1+|\eta|)^N(1+|z|)^{-N}
            e^{(r/2)|\operatorname{Im}z|}.
\]
Hence Theorem~\ref{PWSDK} yields
\[
 g_\eta:=\Fk^{-1}\bigl(F(\cdot-\eta)\bigr)
 \in C_c^\infty(B_r).
\]
The support bound is independent of $\eta$.

It remains to bound the weighted spectral norm relative to
$\|g_\eta\|_{2,k}$ uniformly in $\eta$.
Put $v_\eta=\mu_k(B(\eta,1))$.
By \eqref{EWBOUND},
\[
 w_k(\eta+\xi)
 \le C v_\eta\langle\xi\rangle^{2\gamma_k},
 \qquad
 \mu_k(B(\eta,\rho))
 \ge c\rho^{d_k}v_\eta.
\]
The constants are independent of $\eta$.
Since $2m>d_k=d+2\gamma_k$, it follows that
\[
 \begin{aligned}
 \bigl\|W(|\cdot-\eta|)\Fk g_\eta\bigr\|_{2,k}^2
 &=\int_{\R^d}W(|\xi|)^2|F(\xi)|^2
                      w_k(\eta+\xi)\,d\xi\\
 &\le CA^2v_\eta
          \int_{\R^d}\langle\xi\rangle^{-2m+2\gamma_k}\,d\xi
 \le C_1v_\eta.
 \end{aligned}
\]
On the other hand, Plancherel and the lower bound for $F$ give
\[
 \|g_\eta\|_{2,k}^2
 =\int_{\R^d}|F(\xi)|^2w_k(\eta+\xi)\,d\xi
 \ge\frac14\mu_k(B(\eta,\rho))
 \ge c_1v_\eta.
\]
Here $C_1,c_1>0$ are independent of $\eta$.
Thus
\[
 f_\eta=\frac{g_\eta}{\|g_\eta\|_{2,k}},
 \qquad
 C_W=\max\{1,(C_1/c_1)^{1/2}\}
\]
prove the first assertion.

For prescribed $\gamma,L$, \eqref{EWBOUND} and the ball
inscribed in each cube imply
\[
 \inf_x\mu_k(x+Q_L)\ge cL^{d_k},
 \qquad
 \mu_k(B_r)=\mu_k(B_1)r^{d_k}.
\]
Choose $r>0$ so small that
$\mu_k(B_r)\le(1-\gamma)cL^{d_k}$ and set
$E=\R^d\setminus\overline B_r$.
Then $E$ is $(\gamma,L)$-thick for $\mu_k$, and every
$f_\eta$ constructed above vanishes on $E$.
This proves failure of observability at every spectral center
with the same $E$ and $C_W$.
Sufficiency follows from Theorem~\ref{T91}.
\end{proof}

\section{The fractal uncertainty principle }
\label{sec:fractal-uncertainty}

In this section we follow the ideas of \cite{AC} and \cite{BD}, and, together with the material obtained in Section~3, we complete the proof of Theorem~\ref{thm:fractal-uncertainty}.
 
Section~\ref{sec:BM} constructs uniformly centered damping multipliers
for line-porous spectral sets. Section~\ref{sec:transference} converts
them into observation on spatial holes. Section~\ref{subsec:iteration}
constructs positive cutoffs and proves Theorem~\ref{thm:fractal-uncertainty}
by iteration over scales.
 
\subsection{Dunkl multipliers and uniformly centered damping}
\label{sec:BM}
We transfer the Beurling--Malliavin multiplier of Theorem~\ref{thm:cohen-classical}
through the spherical Paley--Wiener theorem. The construction acts on
entire multipliers and does not use a modulation identity for $\Fk$.

\begin{lemma}\label{lem:PW-transfer}
Suppose $f$ is entire and $|f(z)|\le A e^{a|\operatorname{Im}z|}$.
For every $\varepsilon>0$ there is an entire $q_\varepsilon$ with
$q_\varepsilon(0)=1$, $|q_\varepsilon(x)|\le1$ for $x\in\R^d$, such
that $f_\varepsilon=f q_\varepsilon$ satisfies \eqref{eq:PW} with radius $a+\varepsilon$.
Consequently $\psi=\Fk^{-1}f_\varepsilon\in C_c^\infty(\overline B_{a+\varepsilon})$.
Real translation $f_\varepsilon(\xi-\eta)$ preserves that support radius for its
inverse Dunkl transform.
\end{lemma}
\begin{proof}
Choose $\varphi\in C_c^\infty(B_\varepsilon)$ nonnegative with
$\int\varphi=1$ and set
$q_\varepsilon(z)=\int\varphi(x)e^{-i\langle x,z\rangle}\,dx$.
Classical Paley--Wiener estimates give
$|q_\varepsilon(z)|\le C_N(1+|z|)^{-N}e^{\varepsilon|\operatorname{Im}z|}$.
Multiplying by the bound on $f$ and applying Theorem~\ref{PWSDK} proves
the assertion. For real $\eta$ use
$1+|z|\le(1+|\eta|)(1+|z-\eta|)$.
The constants in the entire estimates may depend on $\eta$; the support
radius does not. On a fixed small real ball, $|q_\varepsilon|\ge1/2$.
\end{proof}

\begin{theorem}\label{thm:DBM}
Assume \eqref{eq:Cohen-hyp}. For every $a>0$ there are
$t,b,c>0$ and $\psi\in C_c^\infty(B_a)$ such that $f=\Fk\psi$ obeys
\begin{equation}\label{eq:DBM}
 |f(0)|\ge b,\qquad |f(\xi)|\le e^{t\omega(\xi)},\qquad
 |f(\xi)|\le e^{-c\langle\xi\rangle^{1/2}}.
\end{equation}
The constants depend only on $a,d,C_{\rm reg},C_{\rm gr}$, and
$f$ has a uniform positive lower bound on a sufficiently small ball
about zero. No $G$-invariance of $f$ or $\omega$ is needed.
\end{theorem}
\begin{proof}
Choose $t>0$ so that the classical theorem applied to $t\omega$ has
bandwidth less than $a/2$. Its output $f_0$ has an entire extension with
$|f_0(z)|\le A e^{(a/2)|\operatorname{Im}z|}$: the compactly supported
ordinary Fourier transform of $f_0$ is integrable by Cauchy--Schwarz.
Apply Lemma~\ref{lem:PW-transfer} to $f_0$ with smoothing radius $\varepsilon<a/2$, and set $f=f_0q_\varepsilon$.
On the real axis the smoothing factor has modulus at most one, and at
zero it equals one, proving \eqref{eq:DBM}.
The ordinary Fourier bandwidth and the global decay of $f_0$ give a
uniform bound on $\nabla f_0$ by Cauchy--Schwarz on the spectral ball.
The fixed smoothing factor has bounded derivative, so continuity is
quantitative and yields the stated lower bound near zero.
\end{proof}

We next construct damping multipliers at arbitrary frequency centers,
following \cite[Proposition~3.5]{AC} and the transfer argument above.

\begin{definition}\label{def:damping}
For $a,b,\rho>0$, an exponent $N>d/2$, and a weight $W$, a damping
multiplier for $Z$ centered at $\eta$ is $m_\eta=\Fk\psi_\eta$ with
$\psi_\eta\in C_c^\infty(B_a)$ such that
\begin{equation}\label{eq:damping}
 \begin{aligned}
 |m_\eta(\xi)|&\ge b &&(|\xi-\eta|\le\rho),\\
 |m_\eta(\xi)|&\le\langle\xi-\eta\rangle^{-N} &&(\xi\in\R^d),\\
 |m_\eta(\xi)|&\le W(|\xi-\eta|)^{-1} &&(\xi\in Z).
 \end{aligned}
\end{equation}
\end{definition}
Separating the lower bound $b$ from the covering radius $\rho$ avoids
ambiguity in squared norm estimates.

\begin{proposition}
\label{prop:damping}
Fix $a>0$, $N>d/2$, $\mu\ge1$, $A_*\ge1$, and $0<\nu\le1/4$.
There exist $\alpha\in(0,1)$ and $b,\rho,c>0$ with the following property.
If $H\ge1$, $Z\subset\overline B(0,A_*H)$ is $\nu$-porous on lines from scales $\mu$ to $H$, then for every $\eta\in\R^d$ there is a damping multiplier
as in \eqref{eq:damping}, with $W=W_{\alpha,c}$.
All constants are independent of $H,Z,\eta$. They may depend on
$a,N,\mu,A_*,\nu,d$. If $H<\mu$, the bounded-diameter case is included
without a porosity requirement.
\end{proposition}
\begin{proof}
We give the details needed for uniformity in the center and the starting
scale, supplementing \cite[Proposition 3.5 and Section 3.5]{AC}.
Write $Z_\eta=Z-\eta$, and first suppose $|\eta|\le2A_*H$ and
$H\ge\mu$. Then $Z_\eta\subset\overline B(0,3A_*H)$.
Let $0<\zeta\le1$ be the exponent in \eqref{eq:line-length} for porosity $\nu/2$, set
$s=1/5$, and choose
\begin{equation}\label{eq:damping-exponents}
 \alpha=1-s\zeta/2.
\end{equation}
Thus $3s<\alpha<1$ and $\alpha+s\zeta>1$.
For $j\ge j_0$, cover $\{2^j\le|x|\le2^{j+1}\}$ by boundedly
overlapping cubes of width $\ell_j=2^j j^{-s}$, with smooth bumps
$\chi_Q$ supported in the cubes and equal to one on their half-cubes.
The half-cubes cover the annulus. Select cubes meeting $Z_\eta$ in
that annulus and define
\[
 \omega_j(x)=-2^j j^{-\alpha}\sum_{Q\ {\mathrm{selected}}}\chi_Q(x),
 \qquad \omega= C_0\sum_{j\ge j_0}\omega_j.
\]
Choose $j_0$, depending only on the fixed parameters, so that
$\ell_j\ge\mu$ and all the following scale comparisons hold for
$j\ge j_0$. Empty families contribute zero. Each bump family is
supported in $\{2^{j-1}\le|x|\le2^{j+2}\}$, and
\begin{equation}\label{eq:omega-derivative}
 |\partial^\beta\omega_j(x)|\le
 C_\beta 2^{j(1-|\beta|)}j^{s|\beta|-\alpha}
       \mathbf1_{Y_j}(x),\qquad |\beta|\le3,
\end{equation}
where $Y_j\subset Z_\eta+B(0,C_d\ell_j)$.
Finite overlap of annuli proves the first bound in
\eqref{eq:Cohen-hyp}.

For occupied annuli, $2^j\le C_{A_*}H$. By \eqref{eq:porosity-thickening}, $Y_j$ is $(\nu/2)$-porous on lines
from scale $2C_d\ell_j/\nu$ and up to $H$.
Each line meets the supporting annulus in a set contained in an interval
of length at most $C2^j$. Cover that interval by $O_{A_*,d}(1)$ segments,
each of length $R_j=\tfrac12\min(2^j,H)$. Since $2^j\le C_{A_*}H$ and
$\ell_j=2^jj^{-s}$, enlarging $j_0$ ensures $2C_d\ell_j/\nu<R_j<H$.
Applying \eqref{eq:line-length} on these segments gives
\begin{equation}\label{eq:omega-line}
 |Y_j\cap\ell|_1\le C2^j j^{-s\zeta}
 \quad\text{for every line }\ell.
\end{equation}
Hence, for $r\asymp2^j$,
\[
 G^*(r)\le C r(\log(e+r))^{-\alpha-s\zeta}.
\]
This proves the second bound in \eqref{eq:Cohen-hyp}, uniformly in $H$.
After fixing the absolute factor $C_0$, the annular lower bound also gives
\begin{equation}\label{eq:omega-on-set}
 \omega(x)\le-c_0|x|/(\log(2+|x|))^\alpha+C_1\quad(x\in Z_\eta),
\end{equation}
where $C_1$ accounts for the fixed ball $|x|\le2^{j_0+2}$.
Its dependence on $\mu$ is retained.

If $|\eta|>2A_*H$, the set $Z_\eta$ lies between the radii
$|\eta|/2$ and $3|\eta|/2$. Only a fixed number of dyadic annuli meet
it. On those annuli use radial bumps of width comparable to $2^j$
and amplitude $-C_0 2^j j^{-\alpha}$. Their derivative bounds satisfy
\eqref{eq:Cohen-hyp}; each annulus contributes at most $Cj^{-\alpha}$
to its growth integral. Since the number of annuli is bounded, the
constants are independent of $\eta$. These bumps give
\eqref{eq:omega-on-set}; annuli below $j_0$ can again be absorbed into
$C_1$. If $H<\mu$ and the center is near, $Z_\eta$ lies in a fixed
ball of radius $3A_*\mu$, and \eqref{eq:omega-on-set} holds even with
$\omega=0$ and a fixed $C_1$. The far-center construction is unchanged.

Apply Theorem~\ref{thm:DBM} to the resulting weight, and denote its multiplier by $m$. Decrease the
multiplier by the fixed factor $e^{-tC_1}$ and, if necessary, by a
further fixed factor so that its global decay is bounded by
$\langle x\rangle^{-N}$. Equation~\eqref{eq:omega-on-set} then gives
the $W_{\alpha,c}^{-1}$ bound for $c=tc_0$ (or a smaller positive $c$).
All these reductions preserve a uniform positive lower bound on a
fixed small ball. Finally put $m_\eta(\xi)=m(\xi-\eta)$ and use
Lemma~\ref{lem:PW-transfer} to obtain its compactly supported inverse
Dunkl transform. This frequency translation is an operation on an
entire multiplier; it does not assert a modulation identity for $\Fk$.
\end{proof}

\begin{remark}\label{rem:damping-dependence}
The damping constants depend on the starting scale $\mu$. In
Theorem~\ref{thm:fractal-uncertainty}, $\nu,d,G$ determine $\mu$
independently of $h$.
\end{remark}

\subsection{From damping multipliers to observability on holes}
\label{sec:transference}
\begin{theorem}\label{thm:transference}
Let $W$ be a unique continuation weight. Assume that $Z$ has damping multipliers at every center with the same parameters
$a,b,\rho,N$ in Definition~\ref{def:damping}, where $N>d/2$.
Suppose $E,S$ are measurable, $E$ is $(\gamma,L)$-thick for $\mu_k$,
and $GE+B_a\subset S$. Then
\begin{equation}\label{eq:holes-observation}
 \supp\Fk f\subset Z\quad\Longrightarrow\quad
 \|f\|_{2,k}\le C_{\rm obs}\|\mathbf1_Sf\|_{2,k}.
\end{equation}
The constant depends only on $W,a,b,\rho,N,\gamma,L,d$, with the fixed
Dunkl structure suppressed. The orbit buffer $GE+B_a\subset S$
ensures convolution localization; no invariance of $E$ or $S$ is required.
\end{theorem}
\begin{proof}
Following the transference scheme of \cite[Theorem 5.2]{BJ} and
\cite[Theorem 3.4]{AC}, take
\begin{equation}\label{eq:lattice}
 \Lambda=(\rho/\sqrt d)\mathbb Z^d,\qquad
 P_\eta f=\Fk^{-1}(m_\eta\Fk f)=f\star_k\psi_\eta.
\end{equation}
Every point is within $\rho/2$ of a lattice point. Pointwise lattice
summation, with respect to the Euclidean lattice rather than the Dunkl
measure, gives
\begin{equation}\label{eq:lattice-bounds}
 b^2\le\sum_{\eta\in\Lambda}|m_\eta(\xi)|^2\le C_\Lambda,
 \qquad
 K_W:=\sup_\xi\sum_{\eta\in\Lambda}W(|\xi-\eta|)^{-1}<\infty.
\end{equation}
The upper bounds follow from $2N>d$ and superpolynomial growth of $W$.
Tonelli and Plancherel therefore yield
\[
 b^2\|f\|_{2,k}^2\le\sum_\eta\|P_\eta f\|_{2,k}^2
 \le C_\Lambda\|f\|_{2,k}^2.
\]
Call $\eta$ good if
$\|W(|\cdot-\eta|)^{1/2}\Fk(P_\eta f)\|_{2,k}
\le A\|P_\eta f\|_{2,k}$, and bad otherwise. Since $\Fk f$ is
supported in $Z$, \eqref{eq:damping} implies
\begin{equation}\label{eq:bad-centers}
 \sum_{\eta\ {\mathrm{bad}}}\|P_\eta f\|_{2,k}^2
 \le A^{-2}\sum_\eta\int_Z W(|\xi-\eta|)|m_\eta(\xi)|^2
                                  |\Fk f(\xi)|^2\,d\mu_k(\xi)
 \le A^{-2}K_W\|f\|_{2,k}^2.
\end{equation}
Choose $A\ge1$ so that $A^{-2}K_W\le b^2/2$.
Since \(W\) is a unique continuation weight, \(W^{1/2}\) is one as well. For each good center, Theorem~\ref{T91} applies to the thick set $E$
with weight $W^{1/2}$. Lemma~\ref{lem:localization}, using
$GE+B_a\subset S$, then gives
\[
 \|P_\eta f\|_{2,k}
 \le C\|\mathbf1_EP_\eta f\|_{2,k}
 =C\|\mathbf1_EP_\eta(\mathbf1_Sf)\|_{2,k}
 \le C\|P_\eta(\mathbf1_Sf)\|_{2,k}.
\]
Consequently
\[
 \frac{b^2}{2}\|f\|_{2,k}^2
 \le\sum_{\eta\ {\mathrm{good}}}\|P_\eta f\|_{2,k}^2
 \le C^2C_\Lambda\|\mathbf1_Sf\|_{2,k}^2.
\]
This proves the assertion, including the squares in the lower bound
and the factor $A^{-2}$ in the bad-center estimate.
\end{proof}

\subsection{Proof of Theorem~\ref{thm:fractal-uncertainty} by iteration}
\label{subsec:iteration}
We first construct positive band-limited cutoffs. Combined with
Theorem~\ref{thm:transference}, they give a fixed contraction at
each admissible scale and hence the power gain in
Theorem~\ref{thm:fractal-uncertainty}.
\begin{lemma}\label{lem:cutoff}
There is a nonnegative radial $\phi\in\mathcal S(\R^d)$ with
$\supp\Fk\phi\subset\overline B_1$ and $\Fk\phi(0)=1$.
Put $\phi_\varepsilon(x)=\varepsilon^{-d_k}\phi(x/\varepsilon)$.
For every bounded $G$-invariant set $X$, $0<\varepsilon<\delta$, and
$A=X+B_\delta$, define
\[
 \Psi=\phi_\varepsilon\star_k\mathbf1_A.
\]
Then $\Psi\in\mathcal S$ is $G$-invariant,
$0\le\Psi\le1$, and $\supp\Fk\Psi\subset\overline B_{\varepsilon^{-1}}$.
For each $M>0$,
\begin{equation}\label{eq:cutoff-estimates}
 \begin{aligned}
 \Psi(x)&\ge1-C_M(\varepsilon/\delta)^M &&(x\in X),\\
 \Psi(x)&\le C_M(\varepsilon/r)^M
       &&(\dist(x,X)\ge\delta+r, r>0).
 \end{aligned}
\end{equation}
\end{lemma}
\begin{proof}
Choose a nonzero real radial $\chi\in C_c^\infty(B_{1/4})$, with
$\chi\ge0$. Its ordinary inverse Fourier transform $q$ is real and
radial. The function $q^2$ is nonnegative, Schwartz, and entire of
exponential type at most $1/2$. Theorem~\ref{PWSDK}, also applied to
the inverse transform, gives $\supp\Fk(q^2)\subset\overline B_{1/2}$.
A positive scalar multiple of $q^2$ has $c_k\int\phi\,d\mu_k=1$.

Let $K_\varepsilon(x,y)=c_k\tau_x\phi_\varepsilon(-y)$.
By \eqref{eq:radial-translation}, $K_\varepsilon\ge0$ and
$\int K_\varepsilon(x,y)\,d\mu_k(y)=1$.
For the radial tail $q_r(z)=\phi_\varepsilon(z)\mathbf1_{\{|z|\ge r\}}$,
the same representation gives
$\tau_x\phi_\varepsilon(-y)=\tau_xq_r(-y)$ whenever $d_G(x,y)>r$.
For positive multiplicities, \eqref{eq:translated-decay} also gives decay
of this kernel away from $Gx$. For all $k\ge0$, radial $L_k^1$
contractivity and scaling give
\begin{equation}\label{eq:radial-tail}
 \int_{d_G(x,y)>r}K_\varepsilon(x,y)\,d\mu_k(y)
 \le c_k\|q_r\|_{1,k}
 \le C_M(\varepsilon/r)^M.
\end{equation}
The assertion for nonsmooth radial tails follows by positive
approximation or the integral representation itself.

Since $A$ is $G$-invariant, if $x\in X$ and $y\notin A$ then
$d_G(x,y)\ge\delta$. If $\dist(x,X)\ge\delta+r$ and $y\in A$,
then $d_G(x,y)\ge r$. Integrating over $A$ and its complement proves
\eqref{eq:cutoff-estimates}; endpoint level sets can be handled by
limits. Covariance gives $G$-invariance of $\Psi$.
Finally $\Fk\Psi=(\Fk\phi_\varepsilon)(\Fk\mathbf1_A)$ is smooth
and compactly supported: $A$ is bounded, so all derivatives of
$\Fk\mathbf1_A$ exist. This proves both the Schwartz assertion and
the spectral support bound.
\end{proof}

\begin{proof}[Proof of Theorem~\ref{thm:fractal-uncertainty}]
\noindent\textbf{Step 1: Symmetrization and observation at each scale.}

Replace $\nu$ by $\min(\nu,1/4)$ and set
\begin{equation}\label{eq:symmetrization}
 X=G\mathbf X,\quad Y=G\mathbf Y,\quad
 \vartheta=(\nu/4)^{|G|},\quad a_0=\vartheta^{-1}.
\end{equation}
Taking closures, if needed, preserves all the open holes and only
enlarges the spatial set and the allowed spectral set. We may therefore
work with compact $X,Y$.
Lemma~\ref{lem:porosity-union} applied to the isometric images shows that $X$ is
$\vartheta$-porous on balls from scales $a_0h$ to $1$ and $Y$ is
$\vartheta$-porous on lines from scales $a_0$ to $h^{-1}$.
Also $X\subset\overline B_{\sqrt d}$ and $Y\subset\overline B_{\sqrt d/h}$.
It suffices to estimate $\mathbf1_Xf$, since $\mathbf X\subset X$,
and the original spectral support is contained in $Y$. We retain these
auxiliary symmetrizations for the positive cutoffs and spectral support
bookkeeping below, rather than as a hypothesis of the observation theorem.

Put $a=\vartheta/16$. Whenever $2a_0h\le r\le1$, let $Q$ range
over all axis-parallel cubes of side $r$. Apply porosity to the
central ball of each $Q$, of diameter $r/2$, and choose a point $x_Q$
with distance at least $\vartheta r/2$ from $X$. The choices need not
be equivariant under $G$. Define
\begin{equation}\label{eq:hole-sets}
 S_r=\{x:\dist(x,X)>ar\},\qquad
 E_r=\bigcup_Q B(x_Q,\vartheta r/8).
\end{equation}
Each chosen ball lies in its cube and in
$\{x:\dist(x,X)>2ar\}$. The union $E_r$ is open, and every cube
of side $r$ contains one of these balls. Thus $r^{-1}E_r$ has fixed
Lebesgue thickness at scale one, and Proposition~\ref{TEV} gives
fixed weighted thickness. Since $X$ is invariant, every point of
$GE_r$ also has distance greater than $2ar$ from $X$. Consequently
\begin{equation}\label{eq:hole-margin}
 G(r^{-1}E_r)+B_{a/2}\subset r^{-1}S_r.
\end{equation}
This is the orbit buffer required in Theorem~\ref{thm:transference};
$E_r$ itself need not be invariant.

For the rescaled spectral set
\begin{equation}\label{eq:rescaled-spectrum}
 Z_r=rY+\overline B_2,\qquad H=r/h,
\end{equation}
\eqref{eq:porosity-thickening} gives porosity with parameter $\vartheta/2$ above the fixed scale
$\mu=4/\vartheta$ and up to $H$. Indeed the starting scale for $rY$
is $ra_0\le a_0$, and thickening by radius two costs at most
$4/\vartheta$. Also $H\ge1$ and
$Z_r\subset\overline B_{(\sqrt d+2)H}$.
Use Proposition~\ref{prop:damping} with support radius $a/2$ and,
for instance, $N=d+1$. It supplies a single weight $W$ and fixed
damping parameters for every $Z_r$ and every frequency center.
Theorem~\ref{thm:transference}, \eqref{eq:hole-margin}, and unitary
rescaling therefore give a constant $c_{\rm obs}\in(0,1]$ such that
\begin{equation}\label{eq:uniform-scale-observation}
 \supp\Fk f\subset Y+\overline B_{2/r}
 \quad\Longrightarrow\quad
 \|\mathbf1_{S_r}f\|_{2,k}\ge c_{\rm obs}\|f\|_{2,k}.
\end{equation}
This constant is independent of $r,h$ and of the large scale-separation
parameter chosen next.

\medskip
\noindent\textbf{Step 2: Cutoffs, spectral support, and contraction.}
Choose $L=2^T\ge2$ later, put $r_j=L^{-j}$, and use
Lemma~\ref{lem:cutoff} with
\begin{equation}\label{eq:iteration-cutoffs}
 \begin{gathered}
 \varepsilon_j=r_j/L,\qquad\delta_j=r_j/\sqrt L,\qquad
 \Psi_j=\phi_{\varepsilon_j}\star_k\mathbf1_{X+B_{\delta_j}},\\
 f_m=\left(\prod_{j=0}^{m-1}\Psi_j\right)f,\qquad f_0=f.
 \end{gathered}
\end{equation}
The functions $\Psi_j$ are invariant, satisfy $0\le\Psi_j\le1$, and
have spectral support in $\overline B_{L/r_j}$.
Since $Y$ is invariant, Lemma~\ref{lem:localization} inductively gives
\begin{equation}\label{eq:iteration-support}
 \supp\Fk f_m\subset Y+\overline B_{L+L^2+\cdots+L^m}
 \subset Y+\overline B_{2/r_m}\quad(m\ge1).
\end{equation}
At $m=0$ this is also true. The geometric sum involves only the
previous cutoffs. For $\widetilde f_m(x)=r_m^{d_k/2}f_m(r_mx)$,
the dilation formula in \eqref{eq:logvinenko-sereda-scaling} yields
$\supp\Fk \widetilde f_m\subset r_mY+\overline B_2$.
This is the rescaling direction required in \eqref{eq:rescaled-spectrum}.

Take $M=4$ in Lemma~\ref{lem:cutoff} and write
$\varepsilon_L=C_4L^{-2}$. Then
\begin{equation}\label{eq:iteration-bounds}
 \Psi_j\ge1-\varepsilon_L\quad\hbox{on }X.
\end{equation}
For $L$ large enough that $L^{-1/2}\le a/2$, the same lemma shows
$\Psi_j\le C_4(2/(aL))^4$ on $S_{r_j}$.
Increase $L$ so this is at most $1/2$. It follows from
\eqref{eq:uniform-scale-observation} that, at every admissible scale,
\begin{align}
 \|f_{m+1}\|_{2,k}^2
 &\le\|f_m\|_{2,k}^2-\frac34\|\mathbf1_{S_{r_m}}f_m\|_{2,k}^2
 \le q^2\|f_m\|_{2,k}^2,
 \qquad q=\sqrt{1-\tfrac34c_{\rm obs}^2}<1.
 \label{eq:contraction}
\end{align}
The number $q$ was fixed before $L$. We can therefore increase $L$
once more so that $\varepsilon_L<(1-q)/2$, and then fix it permanently.

\medskip
\noindent\textbf{Step 3: Iteration and the power of $h$.}
Suppose first that $2a_0h\le1$, and let
\[
 n=1+\left\lfloor\log_L\frac1{2a_0h}\right\rfloor.
\]
Exactly the scales $r_0,\ldots,r_{n-1}$ used here satisfy
$r_j\ge2a_0h$. Iterating \eqref{eq:contraction} and comparing with
\eqref{eq:iteration-bounds} gives
\[
 (1-\varepsilon_L)^n\|\mathbf1_Xf\|_{2,k}
 \le\|f_n\|_{2,k}\le q^n\|f\|_{2,k}.
\]
Put
\begin{equation}\label{eq:beta}
 \beta=-\frac{\log(q/(1-\varepsilon_L))}{\log L}>0.
\end{equation}
Since $n\ge\log_L(1/(2a_0h))$,
$\|\mathbf1_Xf\|_{2,k}\le(2a_0)^\beta h^\beta\|f\|_{2,k}$.
For $2a_0h>1$, the trivial inequality
$\|\mathbf1_{\mathbf X}f\|_{2,k}\le\|f\|_{2,k}$ gives the same
conclusion after enlarging $C$. This proves \eqref{eq:fractal-uncertainty}.
\end{proof}

\end{document}